\documentclass[journal]{IEEEtran}

\usepackage{mathcomSTEv4}
\usepackage[acronym]{glossaries}
\newacronym{dlmd}{\textsc{{DLMD-DiffEx}}}{Decentralized Lazy Mirror Descent with Differential Exchanges}
\usepackage{cite}
\usepackage{amsmath,amssymb,amsfonts}
\usepackage{enumerate}
\usepackage{bbm}
\usepackage{tikz}
\usepackage{mathabx}
\usepackage{graphicx}
\usetikzlibrary{positioning}
\usetikzlibrary{arrows.meta}
\usetikzlibrary{decorations.markings}
\usetikzlibrary{decorations.pathreplacing,decorations.markings}
\usetikzlibrary{arrows,calc,decorations.markings,math,arrows.meta}
\usetikzlibrary{shapes.multipart}
\usetikzlibrary{arrows,positioning,shapes,backgrounds,fit} 
\usepackage{algpseudocode}
\usepackage{algorithm}
\usepackage{relsize}
\usepackage{multirow}
\usepackage{array}
\usepackage{hyperref}

\usepackage[textfont=normalfont,singlelinecheck=off,justification=centering,font=footnotesize]{caption}

\usepackage[margin=30pt,textfont=normalfont,singlelinecheck=off,justification=raggedright,font=footnotesize]{subcaption}
\usepackage{dsfont}

\usepackage{environ}         
\usepackage{etoolbox}        
\usepackage{graphicx}        

\usepackage{anyfontsize}         

\usepackage{xcolor}
\usepackage{multicol}
\usepackage[utf8]{inputenc}

\usepackage{bbm}
\usepackage{authblk}

\allowdisplaybreaks
\usepackage{tikz,pgfplots}
\usetikzlibrary{arrows}

\newcolumntype{P}[1]{>{\centering\arraybackslash}p{#1}}
\DeclareRobustCommand*{\IEEEauthorrefmark}[1]{%
  \raisebox{0pt}[0pt][0pt]{\textsuperscript{\footnotesize #1}}%
}

\newtheorem{rem}{Remark}
\newtheorem{assumption}{Assumption}

\newtheorem{theorem}{Theorem}
\newtheorem{lemma}{Lemma}

\newcommand{\ntv}{\widetilde{\nv}}

\usepackage{caption}
\usepackage{subcaption}

\begin{document}

\title{\huge{Decentralized optimization over noisy, rate-constrained networks: Achieving consensus by communicating differences}}

\author{\IEEEauthorblockN{Rajarshi~Saha\IEEEauthorrefmark{1},
Stefano~Rini\IEEEauthorrefmark{2},
Milind~Rao\IEEEauthorrefmark{3}, 
Andrea~Goldsmith\IEEEauthorrefmark{1,4}}\\
\IEEEauthorblockA{\IEEEauthorrefmark{1}Stanford University,}
\IEEEauthorblockA{\IEEEauthorrefmark{2}National Yang Ming Chiao Tung University,}
\IEEEauthorblockA{\IEEEauthorrefmark{3}Amazon~Alexa~AI,}
\IEEEauthorblockA{\IEEEauthorrefmark{4}Princeton University}
\thanks{This work is supported by grants ONR N00014-18-1-2191, Intel grant 134571, and MOST grant 110-2221-E-A49-052.}}

\markboth{Preprint}%
{}

\maketitle

\begin{abstract}
In decentralized optimization, multiple nodes in a network collaborate to minimize the sum of their local loss functions. 
The information exchange between nodes required for this task, is often limited by network connectivity.
We consider a setting in which communication between nodes is hindered by both (i) a finite rate-constraint on the signal transmitted by any node, and (ii) additive noise corrupting the signal received by any node.
We propose a novel algorithm for this scenario: \gls{dlmd}, which guarantees convergence of the local estimates to the optimal solution under the given communication constraints.
A salient feature of \gls{dlmd} is the introduction of additional proxy variables that are maintained by the nodes to account for the disagreement in their estimates due to channel noise and rate-constraints.
Convergence to the optimal solution is attained by having nodes iteratively exchange these disagreement terms until consensus is achieved. 
In order to prevent noise accumulation during this exchange,
\gls{dlmd} relies on two sequences: one controlling the power of the transmitted signal, and the other determining the consensus rate.
We provide insights on the design of these two sequences which highlights the interplay between consensus rate and noise amplification.
We investigate the performance of \gls{dlmd} both from a theoretical perspective as well as through numerical evaluations on synthetic data and MNIST\footnote{MATLAB and Python implementations can be found \href{https://github.com/rajarshisaha95/DLMD-DiffEx}{\textcolor{cyan}{here}}.}. 
\end{abstract}

\begin{IEEEkeywords}
Decentralized optimization; Lazy mirror descent; Additive channel noise; Finite data rate constraint.
\end{IEEEkeywords}

\IEEEpeerreviewmaketitle

\section{Introduction}
\label{sec:intro}

\IEEEPARstart{I}{n} today's networked world, enormous amounts of data are being collected by edge devices in a distributed fashion.
It is often preferred to keep the data localized on the edge devices themselves for privacy concerns, fault-resilience, or the sheer communication expense of centrally aggregating the available data. 
As a consequence, there has been a rapid emergence of decentralized computational strategies in which  global learning tasks are accomplished  only through local processing \cite{ieeespm2020}. 
These algorithms often offer convergence guarantees of the local iterates at each node to the optimal solution at rates comparable to centralized approaches, even under the constraint that nodes communicate only with their immediate neighbors. 
In system implementations like wearable sensors \cite{huang2018loadaboost}, IoT for smart agriculture \cite{ayaz2019internet},  or vehicular communications \cite{samarakoon}, edge devices might communicate over a wireless channel, and are thus limited in connectivity and reliability by the radio environment. 
Consequently, in addition to connectivity constraints  introduced by the network topology, communication between nodes is also subject to (i) finite rate-constraints, and (ii) corruption by additive noise.
Although these two constraints are considered separately in the literature, to the best of our knowledge their combination has been overlooked.

In this work, we emphasize that the simultaneous presence of both the constraints gives rise to a challenging situation of noise accumulation at the nodes, which can subsequently cause any decentralized optimization algorithm to diverge.
Ensuring that the simultaneous effect of quantization noise and channel noise accumulation does not affect convergence is not addressed by existing approaches.
For this reason, we propose \gls{dlmd} to perform decentralized optimization in challenging environments of this nature.\footnote{A preliminary version of this work was presented in ICASSP 2021. \cite{saha_icassp_2021}
This work presents a much more detailed analysis of the algorithm and simulations, as well as general sufficient conditions for controlling noise accumulation.
}

\vspace{-2mm}
\subsection{Relevant prior work}
\label{subsec:relevant_prior_work}
The last decade has had a flurry of research on decentralized optimization.
In what follows, we briefly mention the results most closely related to our work in this paper.

Models for decentralized information exchange and computation were first proposed by Tsitsiklis et al. \cite{tsitsiklis_bertsekas_1986_TAC}. 
Building upon their information exchange model, Ned\'ic and Ozdagler \cite{nedic2009TAC} proposed subgradient methods for collaborative optimization of a global objective whose components are distributed across multiple agents. 
This was followed by the work of Duchi et al. \cite{duchi2012TAC} which characterized the $\Ocal\Paren{\loge{K\sqrt{n}}/\sqrt{K}}$ (where $K$ is the number of iterations and $n$ is the number of nodes) convergence rate for the class of primal-dual distributed algorithms for non-smooth Lipschitz continuous convex objective functions.
They also address the dependence of convergence rates on the number of nodes and network topology. 

Subsequent works considered various regularity assumptions on the objective function, and proposed distributed algorithms with different convergence rates under such assumptions. 
For example, Tsianos et al. 
\cite{tsianos2012allerton} obtained a convergence rate of  $\Ocal\Paren{\loge{K\sqrt{n}}/K}$ for strongly convex functions, while Jakovetic et al. \cite{jakovetic2014} showed $\Ocal\Paren{K^{-2}}$ rates for first-order smooth objective functions. 
In our work, we consider the problem of optimizing non-smooth Lipschitz continuous convex objective functions (similar to \cite{duchi2012TAC}) under rate-constrained noisy settings in which nodes have oracle access to (possibly stochastic) subgradients.

The above-mentioned literature often does not account for constraints in communication beyond the nodes' connection topology.
Such constraints play an important role for many scenarios of practical relevance, in which the quality of the communication link is imperfect.
Decentralized algorithms require communication between devices, which is often a major bottleneck to the performance of these algorithms in practice, given wireless and wireline network rate constraints and latency.
Consequently, a lot of research effort has gone into understanding the effects of constraints other than network topology \cite{Koloskova2020Decentralized, Tang_NIPS2018, kovalev_2021_aistats}.
Broadly, two major veins of work can be identified with regard to compressing the messages being exchanged among nodes.
The first line of works introduces \textit{dimensionality constraint} on the size of messages being exchanged to reduce communication demand \cite{Alistarh_NIPS2018, Wangni_NIPS2018, MilindICASSP2018, StefanoAsilomar2019, smith2016cocoa}.
Random sparsification, top-$k$ sparsification, round-robin coordinate updates and related schemes, are common strategies in these works.
The other vein considers a \textit{cardinality constraint} on the messages exchanged and quantizes the information to finite precision \cite{qsgd_alistarh, Wen_2017_TernGrad, pmlr-v80-bernstein18a}.
Some of these works \cite{Alistarh_NIPS2018, qsgd_alistarh} consider parameter server frameworks for distributed optimization, as opposed to fully decentralized; nevertheless we mention them because their information compression schemes are invariably applicable to decentralized networks as well.

For wireless environments, quantization is essential because nodes cannot exchange information with infinite precision over a channel with finite bandwidth, as noted in \cite{mmamiri2020}. 
Some works on decentralized optimization that consider quantization \cite{qsgd_alistarh, nedic_quantization, yuan_20121053, rabbat_quantized_incremental}, have quantizers with infinite dynamic range. 
For such infinite-level quantizers, the quantization error is always bounded.
In practice however, for channels with finite bandwidth, quantizers must always have a finite dynamic range. 
If the input to a finite-range quantizer exceeds allowable bounds, it can get saturated, resulting in the quantization error getting unbounded and subsequent divergence of the algorithm.
Existing works either ignore these boundary effects, or even if some works like \cite{yi2014TCNS, Reisizadeh_2019} do take into account a finite dynamic range, they consider a noiseless environment, as a result of which they are unambiguously able to choose large enough initial dynamic ranges so that the quantizers never get saturated throughout the execution of the algorithm.

The presence of additive noise affecting the information exchange among nodes is not well studied, except for \cite{srivastava2010}. 
Our work considers quantization in the presence of channel noise for which we address unique challenges and propose \gls{dlmd}.

\vspace{-1.5mm}
\subsection{Our contributions}
\label{subsec:our_contributions}

To the best of our knowledge, this work is the first to develop an algorithm for distributed optimization over a network comprised of both (i) a finite data rate constraint over information exchanged between nodes, and (ii) corruption of communications between nodes by channel noise.
The interplay between these two communication constraints presents a novel set of challenges that is not addressed with existing approaches. 
We propose \gls{dlmd}, a variation of lazy mirror descent (a.k.a. Nesterov’s dual averaging algorithm \cite{nesterov_2009_mathprog, duchi2012TAC}) to address this issue.
The salient features of \gls{dlmd} are:
\begin{enumerate}[1.]
    \item For achieving consensus between models across different nodes, \gls{dlmd} exchanges differential updates instead of the whole model.
    At any point of time, each node maintains an estimate(s) of the model(s) of each of its neighbor(s).
    These running estimates are iteratively updated as nodes receive information about how their neighbors' current (updated) state \textit{disagrees} from its own past estimate.

    \item These differential exchanges are quantized to a finite set of values and their amplitudes scaled appropriately before transmission, in order to satisfy the rate-constraints and dictate power requirements.

    \item With these quantized differentials being exchanged over a noisy channel, the contribution of additive channel noise accumulates as the iterations proceed, subsequently causing the algorithm to diverge. 
    We introduce \textbf{confidence} and \textbf{power-control} sequences for \gls{dlmd} to prevent this noise accumulation.
    The confidence sequence dictates the eagerness with which a node relies on the information it receives from its neighbors, which in turn controls the convergence rate of the algorithm. On the other hand, the power-control sequence specifies the aggressiveness with which noise accumulation is actively mitigated through amplitude scaling over successive iterations. 
    \item  We also focus on the design of finite-rate quantizers for quantizing the differential exchanges 
    and show that by choosing the dynamic range of the quantizer appropriately, we can simultaneously ensure that the algorithm converges (albeit at a slower convergence rate) while also ensuring that the quantizers remain unsaturated with probability arbitrarily close to one.

\end{enumerate}

Our main result is the derivation of an upper bound to the expected suboptimality gap after $K$ iterations, that scales as $\Ocal\left(K^{-\frac{1-\gamma}{2}}\right)$, $0 < \gamma \leq 1$, given that the quantizers remain unsaturated with sufficiently high probability.
Note that depending on the choice of $\gamma$, this asymptotic rate is close to $\Ocal\left(K^{-1/2}\right)$, which is characterized to be the optimal rate for the convergence of subgradient-type algorithms for the class of non-smooth convex objective functions \cite{duchi2012TAC}.
As a corollary to these results, we show how the confidence and power control sequences can be jointly designed so as to trade-off between the convergence speed and the power requirements of the algorithm. 

We remark that our information exchange model over networks is quite general and extends beyond wireless settings.
It applies to a number of edge computation scenarios where communication among devices occurs over a pre-determined physical layer with a given symbol rate and symbol error probability.
We do not consider any assumption on the additive noise source, except for finite variance.
For wireline environments, it can be any source of external noise such as the decoding error at the receiver (with or without any channel coding scheme).
It can also include any deliberately introduced disturbances in the presence of an adversary with finite power.

\vspace{-3mm}
\subsection{Notation}
\label{subsec:notations}

Lower case boldface letters (eg. $\zv$) are used for column vectors and uppercase boldface letters (eg. $\Pv$) designate matrices. 
We also adopt the shorthands $[m:n] \triangleq \{m, \ldots, n\}$
and  $[n] \triangleq \{1, \ldots, n\}$. 
Subscripts are used for node indices and parentheses depict iteration index (eg. $\zv_i(k)$ denotes the state of node $i$ at iteration $k$).
$\overline{\Paren{\cdot}}$ denotes the average of corresponding quantities across different nodes and $\widehat{\Paren{\cdot}}$ is used for iteration-averaged quantities. 
The noisy version of a variable is indicated with a tilde, i.e. $\widetilde{\av}$ is the version of $\av$ after being corrupted by noise.
For any arbitrary norm $\norm{\cdot}$, its dual norm is defined as $\norm{\cdot}_* = \sup_{\xv: \norm{\xv} \leq 1} \inprod{\cdot, \xv}$, where $\inprod{\cdot,\cdot}$ denotes the inner product.
We primarily consider Euclidean norm  $\norm{\cdot}_2$ in this work, which is its own dual.
Given a sub-differentiable function $f$, its sub-differential at $\xv$ is denoted by $\partial f(\xv)$.
Scripts are used to denote sets (eg. $\Acal$) and $\abs{\Acal}$ is used to denote its cardinality.

\vspace{-1.5mm}
\subsection{Organization of the paper}
\label{subsec:organization_of_the_paper}
The remainder of the paper is organized as follows. 
Sec. \ref{sec:system_model} formally introduces the problem statement of decentralized optimization along with our communication model over a network. 
Sec. \ref{sec:preliminaries} summarizes the necessary preliminaries of the Distributed Lazy Mirror Descent (DLMD) algorithm from \cite{duchi2012TAC} and explains the challenges that arise in the presence of finite rate-constraint and channel noise. 
Our proposed algorithm, \gls{dlmd} is presented in Sec. \ref{sec:proposed_algorithm}.
The convergence theorems of \gls{dlmd} are given in Sec. \ref{sec:convergence_analysis} and a proof sketch in Sec. \ref{sec:proofs_sketch}.
Simulation results are provided in Sec. \ref{sec:numerical_simulations}.
After some additional remarks in Sec. \ref{sec:additional_remarks}, the paper concludes in Sec. \ref{sec:conclusions}.
The detailed proofs are given in the appendices.

\section{System model}
\label{sec:system_model}

We first formally introduce the problem of decentralized optimization in Sec. \ref{subsec:decentralized_optimization_setting}, followed by our communication model in Sec. \ref{subsec:communication_setting} pertinent for the system implementation of any decentralized optimization algorithm over a noisy, rate-constrained network. 
Then in Sec. \ref{subsec:performance_evaluation}, we introduce the performance evaluation metrics  of any optimization algorithm subject to the constraints of our communication model.

\vspace{-2mm}
\subsection{Decentralized optimization setting}
\label{subsec:decentralized_optimization_setting}

Consider a network of $n$ computation nodes (or edge devices), each of which is privy to a local objective function $f_i:\Real^d \to \Real$.
In distributed optimization, the goal is to solve the following minimization problem:
\begin{equation}
     \label{eq:distributed_optimization_problem_statement}
     \minimize_{\xv \in \Xcal} \frac{1}{n}\sum_{i \in [n]}f_i(\xv),
\end{equation}
where, $\Xcal$ is a closed, convex set and each $f_i$ is taken to be convex and sub-differentiable, but need not be smooth.
Furthermore, we also assume that there exists a finite optimum $\xv^* \in \Xcal$ of \eqref{eq:distributed_optimization_problem_statement} and each $f_i$ is $L$-Lipschitz continuous with respect to some norm $\norm{\cdot}$, that is, $\abs{f_i(\xv) - f_i(\yv)} \leq L \norm{\xv - \yv} \hspace{2mm} \text{for all} \hspace{2mm} \xv, \yv \in \Xcal$.
Our objective is to solve \eqref{eq:distributed_optimization_problem_statement} in a fully decentralized setting, i.e. there is no central server to orchestrate the operation of the different nodes/agents.
For simplicity, we assume that all nodes  have a common clock so that they can exchange information with their neighbors in a synchronous fashion.\footnote{Note that this assumption can be relaxed using techniques from \cite{rabbat2012cdc_delays}.}
In distributed optimization algorithms for solving \eqref{eq:distributed_optimization_problem_statement}, each node $i \in [n]$ maintains its own estimate of the optimal solution of \eqref{eq:distributed_optimization_problem_statement} at iteration $k$, which is denoted as $\xv_i(k)$.
Iterative algorithms as in \cite{duchi2012TAC} ensure that $\xv_i(k) \to \xv^*$ as $k \to \infty$ for all $i \in [n]$.

\begin{rem}
\label{rem:stochastic optimization}
The framework \eqref{eq:distributed_optimization_problem_statement} is quite general.
For example, it encompasses the stochastic optimization problem in machine learning in which node $i \in [n]$ observes data points $\{ \xi_i^j\}_{j \in [m]}$ drawn from some probability distribution $\Pcal$. Suppose $\Xcal$ denotes the parameter space and $\ell\Paren{\xv, \xi}$ denotes the loss incurred with a particular model $\xv \in \Xcal$ for a random data point $\xi \sim \Pcal$. One is then interested in minimizing the \textit{population risk}: $\minimize_{\xv \in \Xcal} L(\xv) \equiv \minimize_{\xv \in \Xcal} \mathbbm{E}_{\xi \sim \Pcal}\Br{\ell\Paren{\xv;\xi}}$.
Since the probability distribution $\Pcal$ is unknown, \textit{empirical risk minimization} resorts to solving the following optimization problem instead:
\[\minimize_{\xv \in \Xcal} \frac{1}{nm}\sum_{i \in [n]}\sum_{j \in [m]}\ell\Paren{\xv;\xi_j^i} = \minimize_{\xv \in \Xcal}\frac{1}{n}\sum_{i \in [n]}f_i(\xv),\] 
where $f_i(\xv) = \frac{1}{m}\sum_{j \in [m]}\ell\Paren{\xv;\xi_j^i}$ is known only to node $i$.
\end{rem}

\vspace{-3mm}
\subsection{Communication setting}
\label{subsec:communication_setting}

\begin{figure}[ht]
	\centering
	{
		\begin{tikzpicture}[
		node distance= 6em and 4em,
		sloped]
		\tikzset{
			box/.style = {
				fill=blue!15,
				shape=rectangle, 
				rounded corners,
				draw=blue!40, 
				align=center,
				text = black,
				font=\fontsize{10}{10}\selectfont},
			dummybox/.style = {
				shape=circle,  
				align=center, 
				minimum size={width("rrrrrrrrrrr")+2pt}},
			arrow/.style={
				color=black,
				draw=blue,
				-latex,
				font=\fontsize{8}{8}\selectfont},
			state/.style={
				fill={blue!10},
				rounded corners,
				draw=blue, thick,
				minimum height=2em,
				text width=2cm,
				inner sep=2pt,
				text centered,
			}
		}
	\node  (node_i) [state] at (0,0) {\textbf{Node $i$} \\ $f_i(\xv)$};
	\node  (node_j) [state]  at (6,0) {\textbf{Node $j$} \\ $f_j(\xv)$};
	\node (plus) [draw, circle] at (3,0) {$+$};
	\draw[->,line width=1pt] (node_j) -- (plus);
	\draw[->,line width=1pt] (plus) -- (node_i);
	\node (noise) at (3,1.2) {{$\nv_{ij}(k)$}};
	\draw[->,line width=1pt] (noise) -- (plus);
	\node (transmit_signal) at (4.25,0.3) {\small{$\sv_{ij}(k)$}};
	\node (received_signal) at (2.1,0.3) {\small{$\rv_{ij}(k)$}};
	\end{tikzpicture}}
	\caption{Information exchange corrupted by additive channel noise in \eqref{eq:transmit_receive_node_ij}.}
	\vspace{-0.5cm}
	\label{fig:two_node_communication_model}
\end{figure}
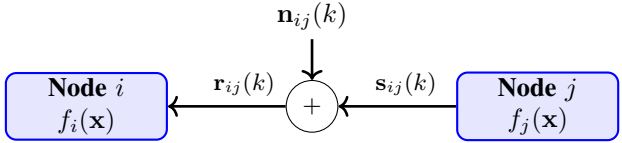

Any decentralized algorithm that solves \eqref{eq:distributed_optimization_problem_statement}, requires communications between the agents in order to guarantee convergence to the optimal solution. 
Similar to prior works, the connectivity of the network is described by an undirected graph $G = (V,E)$, where $V = [n]$ is the set of nodes, each corresponding to a particular agent, and $E \subseteq [n] \times [n]$ denotes the set of edges such that $(i,j) \in E$ if and only if nodes $i$ and $j$ are connected via a link.\footnote{Typically, if the communication latency between two agents and/or the noise corruption are small (for example due to close proximity), then the corresponding nodes of $G$ are considered to be connected by an edge.}
Also, let $\Ncal(j) = \{k \in [n] : (j,k) \in E\}$ be the neighborhood of node $j$.
At iteration $k$, information of the transmitting node $j$ is encoded as $\sv_{ij}(k) \in \Real^d$, transmitted over the channel, and is received at node $i$ (see Fig. \ref{fig:two_node_communication_model}) as:
\begin{equation}
    \label{eq:transmit_receive_node_ij}
    \rv_{ij}(k) = \sv_{ij}(k) + \nv_{ij}(k),
\end{equation}
where, $\nv_{ij}(k)$ is the additive channel noise vector with zero mean $\Expect\Br{\nv_{ij}(k)} = \mathbf{0} \in \Real^d$ and bounded covariance with uncorrelated entries, i.e. $\Expect\Br{\nv_{ij}(k)\nv_{ij}(k)^\top} = \sigma^2\Iv_d$.\footnote{We assume that the links between nodes are non-interfering. In practice, this can be ensured through various multiplexing techniques and is left for future research.}
The \textit{transmit power consumption} for $K$ iterations for transmissions over edge $(i,j) \in E$ is given by: 
\begin{equation}
    \label{eq:power_requirement}
    \Psf(K) = \frac{1}{K}\sum_{k \in [K]}\norm{\sv_{ij}(k)}_2^2.
\end{equation}
Furthermore, since nodes cannot send information with infinite precision and need to exchange quantized values, we assume that $\sv_{ij}(k)$ can only take values in a finite constellation set $\Scal_k \subset \Rbb^d$. 
The \textit{cardinality constraint} on $\Scal_k$ due to quantization is taken to be:
\begin{equation}
    \label{eq:cardinality_constraint}
    \log_2\abs{\Scal_k} \leq \Rsf d,
\end{equation}
where, $\Rsf$ is the \textit{number of bits} per dimension.

\subsection{Performance evaluation}
\label{subsec:performance_evaluation}

Our objective is to design an algorithm that solves the optimization problem in 
Sec. \ref{subsec:decentralized_optimization_setting}
while communication among nodes takes place over links modeled as in Sec. \ref{subsec:communication_setting}.
In this section, we address the question of how to judge the performance of an algorithm under such constraints.
For a given iteration horizon $K$, the \textit{expected suboptimality gap} is defined as:
\begin{equation}
    \label{eq:suboptimality_gap}
    \esf(K) \triangleq \quad  \max_{i \in [n]} \ \Expect [f\Paren{\xvh_i(K)}] - f(\xv^*),
\end{equation}
where, $\xvh_i(K) = \frac{1}{K}\sum_{k \in [K]}\xv_i(k)$ is the iteration-averaged model iterate at node $i$, and $f\Paren{\xv^*}$ is the optimal value.
The expectation in \eqref{eq:suboptimality_gap} is taken over all sources of stochasticity 
(e.g. channel noise) that come into picture.
We assume that, if at any time during the execution of the algorithm up until horizon $K$, constraint \eqref{eq:cardinality_constraint} is violated, then an alarm message is raised and the algorithm terminates.\footnote{In practice, the constraint in \eqref{eq:cardinality_constraint} is satisfied through finite range quantization. For this reason, \eqref{eq:cardinality_constraint} is violated if the quantizer input exceeds its dynamic range, thus leading to saturation.}
We also consider the probability of this event never happening, which is denoted as $P_{\rm S}(K)$ (The subscript $S$ is for \textit{success}).
There exists a trade-off between these two quantities, and we design the algorithm so that $\esf(K) < \epsilon$ for some tolerance level $\epsilon$, and $P_{\rm S}(K) > 1 - \nu$ for some threshold $\nu$.

In other words, for a particular application, given a lower bound on the {success probability}, we want to optimize the {expected suboptimality gap}, i.e.
\begin{equation}
    \label{eq:residual_gap_success_prob}
    \Rcal \lb K, \nu \rb \triangleq \min_{ P_{\rm S}(K)>1-\nu} \esf(K).
\end{equation}

\section{Preliminaries}
\label{sec:preliminaries}
We now briefly review some of the relevant concepts of decentralized lazy mirror descent (a.k.a. distributed dual averaging from \cite{duchi2012TAC}), and highlight the challenges that come into the picture in its implementation over a noisy network with finite rate constrained links.
In Sec. \ref{subsec:Decentralized_Lazy_Mirror_Descent} and \ref{subsec:finite_range_scalar_quantizers}, we respectively review the decentralized version of lazy mirror descent, and present some background on finite-range quantization. 
In Sec. \ref{subsec:implementation_challenges}, we comment on the role of quantization in adapting decentralized lazy mirror descent to the communication channel model that was presented in Sec. \ref{subsec:communication_setting}.

\vspace{-2mm}
\subsection{Decentralized Lazy Mirror Descent (DLMD) algorithm}
\label{subsec:Decentralized_Lazy_Mirror_Descent}

\begin{algorithm}
	\caption{
    Decentralized Lazy Mirror Descent
	}
	\label{algo:DLMD_pseudocode}
	\begin{algorithmic}[1]
		\State {\bf Input:} $K$, $\Pi_{\Xcal}^{\psi}(\cdot)$, $\{\eta(k)\}$
		\State {\bf Initialization:} $\zv_i(0) = \mathbf{0}$, $\xv_i(0) = \mathbf{0}$
		\For {$ k \in [K]$}
		\State \textbf{Consensus step}: $\zv_i'(k+1) = \sum_{j \in [n]}\Pv_{ij}\zv_j(k)$
		\State \textbf{Subgradient step}: $\zv_i(k+1) = \zv_i'(k+1) + \gv_i(k)$
		\State \textbf{Projection step}: $\xv_i(k+1) = \Pi_{\Xcal}^{\psi}\Paren{\zv_i(k+1), \eta(k)}$
		\EndFor
		\State {\bf return}  $\lcb \xv_i(K) \rcb_{i \in [n]} $
	\end{algorithmic}
\end{algorithm}

Consider a doubly stochastic matrix $\Pv \in \Real^{n \times n}$ that respects the structure of the network graph $G = (V,E)$, i.e. for $i \neq j$, $\Pv_{ij} > 0$ iff $(i,j) \in E$, $\Pv_{ij} = 0$ otherwise, and $\sum_{j \in [n]}\Pv_{ij} = \sum_{i \in [n]}\Pv_{ij} = 1$. 
In decentralized lazy mirror descent algorithm, also referred to as distributed dual averaging \cite{duchi2012TAC}, each node $i \in [n]$ maintains a state/model $\zv_i(k)$ at iteration $k$. 
The pseudocode of DLMD \cite{duchi2012TAC} is given in Alg. \ref{algo:DLMD_pseudocode}.
Each iteration of the algorithm consists of three steps: $(i)$ \textbf{Consensus}: Every node $i$ computes a consensus of its neighboring nodes, $(ii)$ \textbf{Subgradient}: Each node $i$ computes a subgradient of its local function $\gv_i(k) \in \partial f_i(\xv_i(k))$, and $(iii)$ \textbf{Projection}: Node $i$ computes a projection of the current state onto the feasible set $\Xcal$ to get the primal iterate $\xv_i(k+1)$.

Here, $\{\eta(k)\}_{k\in \Nbb}$ is a decreasing step size sequence.
$\psi:\Xcal \to \Real$ is a $1$-strongly convex proximal function with respect to some norm $\norm{\cdot}$, i.e. it satisfies
\ean{
\psi(\yv) \geq \psi(\xv) + \inprod{\bigtriangledown \psi(\xv), \yv-\xv} + \f 12\norm{\xv - \yv}^2,
}
for $x,y \in \Xcal$, $\psi \geq 0$ over $\Xcal$ and $\psi(0) = 0$.
Examples include: $(i)$ The \textit{quadratic function} $\psi(\xv) = \f 12\norm{\xv}_2^2$ for $\xv \in \Real^d$ w.r.t. $\ell_2$-norm, and $(ii)$ the \textit{entropy function}: $\psi(\xv) = \sum_{i\in [d]}\xv_i \loge{\xv_i} - \xv_i$ over the probability simplex $\{\xv \hspace{1mm} \vert \hspace{1mm} \xv \succcurlyeq 0, \sum_{i\in [d]}\xv_i = 1\}$ w.r.t. $\ell_1$-norm.
The projection operator $\Pi_{\Xcal}^{\psi}\Paren{\cdot}$ is defined as 
\ean{
\Pi_{\Xcal}^{\psi}\Paren{\zv,\eta} = \argminimize_{\xv \in \Xcal}\left\{\inprod{\zv, \xv} + \frac{1}{\eta}\psi(\xv) \right\}.
}
The role of the proximal function is to ensure that the \textit{primal iterate} $\xv_i(k)$ does not oscillate wildly. 
We refer the reader to  \cite{duchi2012TAC} for a more comprehensive exposition. 
In the following, we consider a modification of decentralized lazy mirror descent in which  steps $(ii)$ and $(iii)$ are substantially unchanged (since they only involve  local computations)  while step $(i)$ is adapted to the communication constraints in Sec. \ref{subsec:communication_setting}.

\subsection{Finite-range scalar quantizers}
\label{subsec:finite_range_scalar_quantizers}

The \textit{consensus step} in Alg. \ref{algo:DLMD_pseudocode} requires nodes to exchange information with each other. 
In practical scenarios, this transmission among nodes requires quantization. 
For the sake of simplicity and ease of implementation, we focus on probabilistic scalar quantizers to illustrate the system design challenges.
Scalar quantizers act on each coordinate of a vector separately.
However, the same concepts are applicable in general for vector quantizers too.

Suppose that a link between any two nodes of the network has limited transmission capabilities, i.e. any message sent over the link is constrained to be of length $\Rsf$-bits per dimension.
Along each dimension, we have $M = 2^\Rsf$ quantization points which we  denote here as $\{u_1, \ldots, u_M\}$, where $u_1 = -U$ and $u_M = +U$.
This means that the dynamic range of the quantizer (denoted by $U$) is finite, and if the input to the quantizer $\vv$ is such that $\norm{\vv}_{\infty} > U$, we say that the quantizer is \textit{saturated}, which is a violation of requirement \eqref{eq:cardinality_constraint}.
The \textit{quantization resolution} is defined as
\[\Delta = \sup_{\vv:\norm{\vv}_{\infty} \leq U}\norm{\vv' - \vv}_{\infty},\] 
where $\vv'$ denotes the quantized version of $\vv$.
Note that for uniform scalar quantizers, $\Delta = 2U/\Paren{2^{\Rsf}-1}$.
We denote our quantization operation by $\Qsf_{\Delta,U}(\cdot)$. 
The probabilistic quantization scheme described next has been employed for distributed consensus problems in \cite{aysal2008TSP, kar_quantization}.
For an input $\vv = [v_1, v_2, \ldots, v_d]^\top \in \Real^d$, the quantizer output is 
\ea{
\vv' = \Qsf_{\Delta, U}\Paren{\vv} = [\qsf(v_1), \ldots, \qsf(v_d)]^\top.
\label{eq:def quant}
}
The operation $\qsf(\cdot)$ is defined for $v \in [u_j, u_{j+1}) \subseteq \Real$, as:
\begin{equation}
    \label{eq:probabilistic_quantizer}
    \qsf(v) = 
    \begin{cases}
           u_j & \text{with probability $r$} \\
           u_{j+1} & \text{with probability $1-r$}, 
    \end{cases}
\end{equation}
where 
\ea{
r = \frac{u_{j+1} - v}{u_{j+1} - u_j} = \frac{u_{j+1} - v}{\Delta}. 
}
It is shown in \cite{aysal2008TSP} that the probabilistic quantization scheme in \eqref{eq:probabilistic_quantizer} is equivalent to a \textit{non-subtractive dithered quantizer} \cite{wannamaker} with uniformly distributed dither sequence.
The following lemma states that as long as $\norm{\vv}_{\infty} \leq U$, the output $\vv' = \Qsf_{\Delta,U}(\vv)$ is an unbiased estimate of $\vv$ and the quantization error has bounded variance.

\begin{lemma}
\label{lem:unbiased_bounded_variance_quantizer}
Suppose $v \in [u_j, u_{j+1}) \subseteq [-U,+U] \subseteq \Real$ and let $v' = \qsf(v)$ be the $\Rsf$-rate quantization of $u$ according to rule \eqref{eq:probabilistic_quantizer}. Then the quantization error $e = v' - v$ has zero mean and bounded variance, that is
$$\Expect\Br{v'} = v, \hspace{2mm} \text{and} \hspace{2mm} \Expect\Br{\Paren{v' - v}^2} \leq \frac{\Delta^2}{4} = \frac{U^2}{\Paren{2^{\Rsf}-1}^2}.$$
\end{lemma}
\begin{IEEEproof}
The probabilistic scalar quantization operation is unbiased since
\ean{\Expect\Br{v'} = u_jr + u_{j+1}\Paren{1-r} = u_{j+1} - \frac{u_{j+1} - v}{\Delta}\Delta = v.
}
To obtain an upper bound on the variance, note that 
\ea{\Expect\Br{\Paren{v' - v}^2} = \Paren{v - u_j}\Paren{u_{j+1} - v}.
\label{eq:prod}
}
Since $v \in [u_j, u_{j+1})$, the product in the RHS of \eqref{eq:prod} is maximized when $v = \frac{u_j + u_{j+1}}{2}$, in which case the maximum attainable value is $\frac{\Delta^2}{4}$.
This completes the proof.
\end{IEEEproof}

\subsection{Implementation challenges}
\label{subsec:implementation_challenges}

Quantization is essential for information exchange in the consensus step.
A naive way to do that would be to directly quantize the states, i.e. node $j$ broadcasts $\Qsf_{\Delta,U}\Paren{\zv_j(k)}$ to its neighbors and any node $i \in \Ncal(j)$  uses the quantized state for computing consensus \cite{nedic_quantization}.
For this strategy to be successful using finite-range quantizers, one requires prior knowledge of model magnitudes. This might not be readily available in many scenarios of interest. 
Alternatively, instead of directly quantizing the model, some works like \cite{yi2014TCNS}, \cite{Reisizadeh_2019} resort to transmitting the quantized difference between the models at iteration $k$ and that at iteration $k-1$.
The intuition behind choosing to do this lies in the fact that when the learning rate is small, nodes need to communicate \textit{less} with their neighbors.
From a theoretical perspective, it is often possible to obtain upper bounds on the $\ell_{\infty}$-norm of the model differentials, which allows us to appropriately design the dynamic range of the quantizers, so that they remain unsaturated.

We now  describe, from a high-level perspective, the challenges that one encounters when quantized differentials are exchanged in the presence of channel noise.
Consider node $i$'s update: Node $i$ requires knowledge of $\zv_j(k)$ for $j \in \Ncal(i)$ to compute consensus, for which both nodes $i$ and $j$ maintain a running estimate $\yv_{ij}(k)$ of $\zv_{j}(k)$.
At iteration $k$, node $j$ sends $\Qsf_{\Delta,U}\Paren{\zv_j(k) - \yv_{ij}(k-1)}$ to node $i$.
In the noiseless case, node $i$ updates its running estimate by adding the quantized differential: $\yv_{ij}(k) = \yv_{ij}(k-1) + \Qsf_{\Delta, U}\Paren{\zv_{j}(k) - \yv_{ij}(k-1)}$, and uses it for consensus.
Note that if $\Qsf_{\Delta,U}\Paren{\cdot}$ had no rate constraint, then $\yv_{ij}(k) = \zv_j(k)$, which is the true state.
Henceforth, $\yv_{ij}(k)$ is referred to as a proxy for $\zv_j(k)$.

The situation becomes more complicated when  differential exchanges between two nodes take place over additive noisy links.
To see this, suppose that the exchange from node $j$ to $i$ is corrupted by an additive noise $\nv_{ij}(k)$.
Node $i$ now receives $\Qsf_{\Delta,U}\Paren{\zv_j(k) - \yv_{ij}(k-1)} + \nv_{ij}(k)$, as a result of which the running estimate of node $j$'s state being maintained at node $i$ is not exactly $\yv_{ij}(k)$, but rather a corrupted version of it, $\ytv_{ij}(k)$.
The contribution from channel noise is the source of discrepancy between the proxies of $\zv_j(k)$ being maintained by nodes $i$ and $j$, i.e. $\ytv_{ij}(k)$ and $\yv_{ij}(k)$ respectively.
Simple algebra (see Appendix \ref{appendix:proof_of_lemma_2}) shows that as iterations proceed, the variance of this discrepancy grows, and is given by
\begin{equation}
    \label{eq:noise_accumulation}
    \ytv_{ij}(k) = \yv_{ij}(k) + \sum_{l\in [k]}\nv_{ij}(l).
\end{equation}
When node $i$ uses these $\ytv_{ij}(k)$'s for consensus, this accumulated noise term can have unbounded variance.
Such noise accumulation can result in the divergence of the naive DLMD algorithm.
The consideration of quantized differential exchanges in the presence of channel noise is unique to our work, and \gls{dlmd} makes use of \textit{diminishing confidences} and \textit{power control} to ensure that the variance of accumulated noise does not grow unbounded. 
In the following section, we propose an algorithm which is capable of ensuring convergence despite the difficulties presented above. 

\subsection{Modeling Assumptions}
Before delving further, we justify our assumptions on the communication setting in Sec. \ref{subsec:communication_setting} as it might appear contradictory that the communication channel has both a rate constraint, typical of digital communication, and is simultaneously corrupted by additive noise, characteristic of analog communication. 
Our model assumes that after quantization, the (digital) quantizer output undergoes a rather general set of transformations before being mapped to a discrete input constellation. 
Post-transmission, the information is demodulated and reconstructed. 
It is assumed that the overall effect of these transformations, whatever they are, can be equivalently expressed as an additive noise of finite variance. 
This is a rather general model, as it allows one to express the overall effect of reconstruction errors arising from the cascade of various operations, including: entropy coding/decoding, channel encoding/decoding with finite block-length, modulation/demodulation, ACK/NACK for varying channel conditions (for example fading).

We emphasize that our distinction of \textbf{quantization noise} and \textbf{additive channel noise} is a fundamental one for distributed optimization algorithms.
The \textbf{quantization noise} can be viewed as a high level abstraction for a variety of noise sources, whose knowledge is present at the transmitting node but not at the receiving node.
This knowledge can be made use of to correct for these sources of error using appropriate feedback mechanisms.
More specifically, for any given communication protocol over a digital link, the noise due to the rate constraint arises from a combination of quantization, modulation, and communication rate.
It can also be due to any compression schemes as in \cite{Koloskova2020Decentralized, Tang_NIPS2018}.

Similarly, the additive error (which we refer to as \textbf{channel noise}) encompasses a broad variety of noise sources, the effects of which are known at the receiving node, but not at the transmitting node, such as demodulation and decoding errors owing to the use of finite block-length codes in communication.
The fundamental question we address in this work is: \textit{How to reconcile the information between the transmitting and receiving nodes when both such noise sources are simultaneously present?}

\section{The \gls{dlmd} algorithm}
\label{sec:proposed_algorithm}

We now introduce \textit{Decentralized Lazy Mirror Descent with Differential Exchanges} (\gls{dlmd}).
The key idea behind \gls{dlmd} is essentially an appropriate modification of the consensus step of DLMD (line 4 in Alg. \ref{algo:DLMD_pseudocode}) in order to avoid the accumulation of noise as in \eqref{eq:noise_accumulation}, for which we introduce (i) \textit{diminishing confidence} and (ii) \textit{power control}.
\begin{enumerate}[(i)]
    \item As iterations proceed, the discrepancy between a node's true state and the estimate being maintained by its neighbor(s) keeps growing.
    Thus, the estimates grow less reliable over time.
    The notion of \textit{confidence} allows for the nodes to progressively decrease the contribution of the neighboring nodes' states in computing consensus and increase reliance on its own state.
    \item In addition to this, \textit{power control} actively increases the transmit power in order to combat the  decreasing effective Signal-to-Noise Ratio (SNR) due to channel noise accumulation. 
    This power control strategy allows \gls{dlmd} to trade-off between the convergence rate and the nodes' transmit power requirement.
\end{enumerate}

\begin{algorithm}
	\caption{
	Decentralized Lazy Mirror Descent with Differential Exchanges (\gls{dlmd}) (Node $i$'s update)
	}
	\label{algo:DLMD-DiffEx_pseudocode}
	\begin{algorithmic}[1]
		\State {\bf Input:} $K$, $\Pi_{\Xcal}^{\psi}(\cdot)$, $\{\eta(k)\}$, $\{\beta(k)\}$, $\{\al(k)\}$
		\State {\bf Initialization:} $\zv_i(1) = \mathbf{0}$, $\xv_i(1) = \mathbf{0}$, $\yv_{ij}(0) = \mathbf{0}$, $\ytv_{ij}(0) = \mathbf{0}$ for all $j \in \Ncal(i)$
		\For {$ k \in [K]$}
		\For { $ j \in \Ncal(i)$ \text{Node} $j$ \text{computes}}
		\If{$\Paren{\norm{\zv_j(k) - \yv_{ij}(k-1)}_{\infty} > U}$}
		\State {\bf return} -- \underline{FAILURE}
		\Else
		\State  Compute quantized state differential:
		\vspace{-2mm}
		\begin{equation*}
		{
		\hspace{16mm}\boldsymbol{\de}_{ij}(k) = \Qsf_{\De,U} \Paren{\zv_j(k) - \yv_{ij}(k-1)}
		}
		\vspace{-2.5mm}
		\end{equation*}
		\EndIf
		\State  Node $j$ updates the proxy for its past state:
		\vspace{-2mm}\ean{\hspace{4mm}\yv_{ij}(k) = \yv_{ij}(k-1) + \boldsymbol{\delta}_{ij}(k)}
		\vspace{-5mm}
	    \State Scale and transmit to node $i \in \Ncal(j)$:
	        \vspace{-2mm}\ean{
			\sv_{ij}(k)= \al(k) \boldsymbol{\de}_{ij}(k)}
		\vspace{-8mm}
		\EndFor
		\State Note $i$ receives $\{\rv_{ij}(k)\}_{j \in \Ncal(i)}$ 
				and decodes:	\vspace{-2mm}\ean{\Tilde{\boldsymbol{\delta}}_{ij}(k) = \alpha(k)^{-1}\rv_{ij}(k)
				}
		\vspace{-6mm}
		\State Node $i$ updates running estimate of $j$'s state: \vspace{-2mm}\ean{\ytv_{ij}(k) = \ytv_{ij}(k-1) + \Tilde{\boldsymbol{\delta}}_{ij}(k)
			}
		\vspace{-4mm}
		\State Compute a (possibly stochastic) subgradient\\ \hspace{26mm}$\gtv_i(k)$ of $f_i$ at $\xv_i(k)$
 		\State \textbf{Consensus/Subgradient step}: Dual update at Node $i$:
 		\vspace{-4mm}
 		\ean{
 		\zv_i\Paren{k + 1} = \Wv_{ii}(k)\zv_i(k) + \hspace{-1.5mm}\sum_{j \in \Ncal(i)}\Wv_{ij}(k)\ytv_{ij}(k) + \gtv_i(k)
 	}
 	 \hspace{15mm} where, $\Wv(k) = \Paren{1 - \beta(k)}\Iv_d + \beta(k)\Pv$
 		\State \textbf{Projection step}: Get primal iterate:
 		\vspace{-2mm}
 		\ean{\xv_i(k+1) = \Pi_{\Xcal}^{\psi}\Paren{\zv_i(k+1), \eta(k)}}
 		\vspace{-7mm}
		\EndFor
		\State {\bf return}  $\lcb \xv_i(K) \rcb_{i \in [n]} $
	\end{algorithmic}
\end{algorithm}

The pseudocode for \gls{dlmd} is given in Alg. \ref{algo:DLMD-DiffEx_pseudocode}.
Each node $i \in [n]$ maintains a \textit{state} or \textit{dual iterate} $\zv_i(k)$ as well as a \textit{primal iterate} $\xv_i(k)$.
Apart from these, it also maintains $\yv_{ji}(k)$ for each of its neighbors $j \in \Ncal(i)$.
$\{\yv_{ji}(k) \hspace{1mm} | \hspace{1mm} j \in \Ncal(i)\}$ is a set of proxies for node $i$'s true state $\zv_i(k)$.
$\yv_{ji}(k)$ is used by node $i$ to compute the differential of its state for node $j$'s update and takes into account the error introduced due to quantization.
Moreover, node $i$ also maintains $\ytv_{ij}(k)$ for $j \in \Ncal(i)$ which are running estimates of its neighboring nodes' states.
The order of subscripts is worth noting: $\yv_{ij}(k)$ is a proxy for $\zv_j(k)$ and is maintained by node $j$ for node $i$, whereas $\ytv_{ij}(k)$ maintained by node $i$, is a noise-corrupted version of $\yv_{ij}(k)$ and acts as an unreliable estimate for $\zv_j(k)$.
Node $i$ uses $\{\ytv_{ij}(k) \hspace{1mm} | \hspace{1mm} j \in \Ncal(i)\}$ in the consensus step of \gls{dlmd} as follows:
\begin{align}
    \label{eq:DLMD-DiffEx_consensus_step}
    &\zv'_i\Paren{k + 1} = \Wv_{ii}(k)\zv_i(k) + \sum_{j \in \Ncal(i)}\hspace{-1.5mm}\Wv_{ij}(k)\ytv_{ij}(k) \nonumber\\
    &= \Paren{1 - \beta(k)\sum_{j \in [n]}\Pv_{ij}}\zv_i(k) + \beta(k)\sum_{j \in [n]}\Pv_{ij}\ytv_{ij}(k),
\end{align}
where  $\Wv(k)=\Paren{1 - \beta(k)}\Iv + \beta(k)\Pv$ is a time-varying sequence of consensus matrices obtained from the doubly stochastic matrix $\Pv$ that describes the network topology, and $\{\beta(k)\}_{k\in \Nbb}$ is the \textit{confidence} sequence. 
The sequence $\{\beta(k)\}$ dictates the rate with which we decay the contribution of the neighboring nodes' states in computing consensus.

Node $i$'s state update at any iteration $k$ requires knowledge of the states of each of its neighbors $j \in \Ncal(i)$.
Node $j \in \Ncal(i)$ computes its \textit{state differential}, $\zv_j(k) - \yv_{ij}(k-1)$ for node $i$ and quantizes it subject to the rate-constraint of the $(i,j)^{th}$ link to get the \textit{quantized differential}, $\boldsymbol{\delta}_{ij}(k) = \Qsf_{\Delta,U}\Paren{\zv_j(k) - \yv_{ij}(k-1)}$.
It then scales it up by $\alpha(k)$ and sends $\sv_{ij}(k) = \alpha(k) \boldsymbol{\delta}_{ij}(k)$ over the noisy link $(i,j) \in E$ to node $i$.
Node $i$ receives $\rv_{ij}(k) = \sv_{ij}(k) + \nv_{ij}(k)$ as per \eqref{eq:transmit_receive_node_ij}.
It then decodes the quantized diferential $\Tilde{\boldsymbol{\delta}}_{ij}(k) = \alpha(k)^{-1}\rv_{ij}(k)$, and updates the running estimate (maintained by itself) of node $j$'s state $\ytv_{ij}(k) = \ytv_{ij}(k-1) + \Tilde{\boldsymbol{\delta}}_{ij}(k)$, which it then uses in computing consensus as in \eqref{eq:DLMD-DiffEx_consensus_step}.
$\{\alpha(k)\}_{k\in \Nbb}$ is the predetermined \textit{power-control} sequence.\footnote{We assume the same sequence $\{\alpha(k)\}$ irrespective of the node just for simplicity. We can have different scaling at different nodes.}
Node $j$ also updates its own proxy for the next iteration, $\yv_{ij}(k) = \yv_{ij}(k-1) + \boldsymbol{\delta}_{ij}(k)$.
The remainder of the algorithm steps involve computing consensus to get $\zv_i'(k+1)$ according to \eqref{eq:DLMD-DiffEx_consensus_step}, evaluating (possibly stochastic) subgradient $\gtv_i(k)$, taking a step in the direction of the subgradient $\zv_i(k+1) = \zv_i'(k + 1) + \gtv_i(k)$, computing a projection to get the primal iterate, $\xv_i(k+1) = \Pi_{\Xcal}^{\psi}\Paren{\zv_i(k+1), \eta(k)}$, and finally updating the time averaged primal iterate.
A summary of the associated variables is given in Table \ref{tab:variable_summary}.

The worst case memory requirement for any individual node $i$ executing \gls{dlmd} is given by $\Ocal\Paren{(2\cdot\text{deg}(i) + 2)d}$, where $\text{deg}(i)$ denotes the degree of node $i$.
The quantity $2\cdot\text{deg}(i)$ comes from the set of proxies $\ytv_{ij}(k)$ and $\yv_{ji}(k)$ for each $j \in \Ncal(i)$.
The additional $2$ comes from node $i$'s own primal and dual iterates $\zv_i(k)$ and $\xv_i(k)$.
Furthermore, apart from the complexity involved in computing consensus, evaluating subgradients, taking a step in the direction of a subgradient, and the projection step (all of which are common to DLMD \cite{duchi2012TAC}), the only additional computation involved is in computing differentials and quantizing them.
For our scalar quantization scheme, this is only $\Ocal(d)$, which is insignificant compared to what the other computations mentioned above involve.
In our work, we consider $\Pv$ to be constant.
Hence, the computational complexity involved with obtaining $\Wv(k)$ is also insignificant.
However, for scenarios where $\Pv$ is varying with $k$, this also needs to be taken into account.
Note that in line 6 of Alg. \ref{algo:DLMD-DiffEx_pseudocode}, when any quantizer saturates, the algorithm does not necessarily ``terminate".
We can, for example, use the current iterate to warm start a fresh execution of our algorithm (after an intermediate model exchange step so that the warm start initialization is same across all the nodes).

\begin{table}[]
\begin{center}
\caption{Summary of variables in \gls{dlmd}}
\renewcommand{\arraystretch}{1.3}
\begin{tabular}{|l|c|}
    \hline
    Primal iterate of node $i$ at iteration $k$ & $\xv_i(k)$ \\
    Dual iterate (state) of node $i$ at iteration $k$ & $\zv_i(k)$ \\
    Proxy for $\zv_j(k)$ maintained by node $j$ for node $i$ & $\yv_{ij}(k)$ \\
    Noisy estimate of node $j$'s state maintained by node $i$ & $\ytv_{ij}(k)$\\
    \hline
    Quantized differential sent from node $j$ to $i$ & $\boldsymbol{\de}_{ij}(k)$ \\
    Noisy quantized differential received by node $i$ from $j$ & $\widetilde{\boldsymbol{\de}}_{ij}(k)$ \\
    \hline
    Signal transmitted by node $j$ to node $i$  &  $\sv_{ij}(k)$  \\
	Signal received by node $i$ from node $j$ &  $\rv_{ij}(k)$\\
	\hline
	(Increasing) Power control sequence & $\al(k)$ \\
	Step-size sequence for projection step & $\eta(k)$ \\
	(Diminishing) Confidence sequence & $\beta(k)$ \\
	\hline
	Quantization resolution & $\Delta$ \\
	Dynamic range of finite rate quantizer & $U$ \\
	\hline
\end{tabular}
\label{tab:variable_summary}
\vspace{-5mm}
\end{center}
\end{table}

\section{Convergence analysis}
\label{sec:convergence_analysis}
In this section, we present our main convergence results, i.e. Thm. \ref{thm:main_convergence_result} and Thm. \ref{thm:quantizer_unsaturation}.
Proof sketches of these theorems are presented in Sec. \ref{sec:proofs_sketch}, while most of the algebra is provided in the appendices.

\subsection{Assumptions}
\label{subsec:assumptions}

We begin with some standard assumptions required for our analysis.
Denote by $\Fcal_k$, the $\sigma$-field containing all information till iteration $k$, that is for all $i \in [n], \gtv_i(1), \ldots, \gtv_i(k)$, $\zv_i(1), \ldots, \zv_i(k+1)$, $\xv_i(1), \ldots, \xv_i(k+1) \hspace{-1.5mm} \in \hspace{-1.5mm} \Fcal_k$. 
Furthermore, for all $i \hspace{-1.5mm} \neq \hspace{-1.5mm} j$, $\yv_{ij}(1), \ldots, \yv_{ij}(k)$, $\ytv_{ij}(1), \ldots, \ytv_{ij}(k)$, $\delv_{ij}(1), \ldots, \delv_{ij}(k)$, $\widetilde{\delv}_{ij}(1), \ldots, \widetilde{\deltv}_{ij}(k) \in \Fcal_k$. 
We consider an oracle model for subgradient evaluations.
Input to the oracle is an evaluation point $\xv \in \Real^d$, and the output is a subgradient $\gv \in \partial f(\xv)$.
Oracle models provide an abstraction and hide the computational complexity involved in subgradient evaluations, which can vary diversely depending on the objective function.

\begin{assumption}
\label{assume:stochastic_subgradient_oracle}
The noisy subgradient $\gtv_i(k)$ is 
{unbiased}, i.e. 
%
%
\ea{
\Expect\Br{\gtv_i(k) \vert \Fcal_{k-1}} = \gv_i(k) \in \partial f_i\Paren{\xv_i(k)},
}
and its magnitude has 
{bounded second moment}, that is 
\ea{
\Expect\Br{\norm{\gtv_i(k)}^2 \vert \Fcal_{k-1}} \leq \Omega^2.
}
\end{assumption}

\begin{assumption}
\label{assume:channel_noise_independence}
The channel noise in the communication model \eqref{eq:transmit_receive_node_ij} is independent of quantization and oracle noise, and is unbiased with $\Var{\Paren{[\nv_{ij}(k)]_s}} = \sigma^2$ where $[\cdot]_s$ denotes the $s^{th}$ coordinate.
\end{assumption}

\subsection{Main results}
\label{subsec:main_results}

One of the notable features of \gls{dlmd} is that it prevents noise accumulation which would occur with naive DLMD in the presence of channel noise, as in \eqref{eq:noise_accumulation}, through an appropriate choice of confidence and power control sequences.

From line 11 of Alg. \ref{algo:DLMD-DiffEx_pseudocode}, we see that  $\sv_{ij}(k)$ scales as $\al(k)$, and hence $\al(k)$ dictates the transmit power requirement of the algorithm as per \eqref{eq:power_requirement}. 
Moreover, \textit{confidence sequence} $\{\beta(k)\}$ is decreasing in $k$, and from \eqref{eq:DLMD-DiffEx_consensus_step}, it is apparent that the rate of decay of the sequence $\beta(k)$ inversely affects the rate at which consensus is achieved, that is the rate at which $\norm{\zvo(k) - \zv_i(k)}_*$ decreases.
Recall that  $\zvo(k) = \frac{1}{n}\sum_{i \in [n]}\zv_i(k)$.
From the aforementioned discussion, we see that the interplay between confidence and power control sequences determines the rate of convergence and the transmit power requirement of \gls{dlmd}.
The following lemma precisely characterizes sufficient conditions on the sequences  $\{\al(k)\}$ and $\{\beta(k)\}$ that prevent channel noise accumulation in \gls{dlmd}.

\begin{lemma}{\bf [Convergence conditions]}
    \label{lem:convergence_rate_power_requirement_trade-off}
    For any choice of $\{\al(k)\}$ and $\{\beta(k)\}$ in \gls{dlmd} (see Alg. \ref{algo:DLMD-DiffEx_pseudocode}), the discrepancy between the noisy estimate of state $j$ maintained by node $i$, i.e. $\widetilde{\yv}_{ij}(k)$ and the proxy for node $j$'s past state maintained by itself, i.e. $\yv_{ij}(k)$ as given by  \eqref{eq:noise_accumulation}, has bounded variance if
    \ea{
    \sum_{l\in [k]}\frac{1}{\al(l)^2} = \Theta\Paren{\frac{1}{\beta(k)^2}}.
    \label{eq:condition al be}
    }
\end{lemma}

Lemma \ref{lem:convergence_rate_power_requirement_trade-off} gives sufficient conditions for general sequences $\{\al(k)\}$ and $\{\beta(k)\}$ to ensure bounded noise accumulation.
The next theorem, Thm. \ref{thm:main_convergence_result} gives an upper bound to the expected suboptimality gap in \eqref{eq:suboptimality_gap} for a particular choice of $\{\al(k)\}$ and $\{\beta(k)\}$ as geometric sequences.
This result holds conditioned on the event that \gls{dlmd} succeeds, i.e. none of the quantizers get saturated. 
Subsequently, Thm. \ref{thm:quantizer_unsaturation} bounds the probability of this happening, and shows that with appropriate choice of design parameters, we can ensure that quantizers remain unsaturated with a sufficiently high probability.

\begin{theorem}
\label{thm:main_convergence_result} [\textbf{Upper bound on expected suboptimality gap}]
Suppose the assumptions in Sec. \ref{subsec:assumptions} hold and the sequence  $\{\xv_i(k)\}_{i \in [n]}$ is generated by \gls{dlmd} as in Alg. \ref{algo:DLMD-DiffEx_pseudocode}, where the \textit{confidence}, \textit{power control}, and \textit{step-size} sequences are chosen as:
\begin{subequations}
\begin{equation}
    \alpha(k) = \sqrt{c_1}k^{\tau/2} \label{power_control_sequence}
\end{equation}
\begin{equation}
    \beta(k) = c_0 k^{-\gamma}
    \label{confidence_sequence}
\end{equation}
\begin{equation}
    \eta(k) = \frac{R\sqrt{1 - \lambda}}{4 \xi k^{(1 + \gamma)/2}},
    \label{step_size_sequence}
\end{equation}
\end{subequations}
where $\xi = \Paren{\Omega^2 + c_0^2\frac{\Delta^2}{4}d + \frac{c_0^2\sigma^2d}{2\gamma c_1}}^{1/2}$ with $\Omega^2$ and $\sigma^2$ given as in Asm. \ref{assume:stochastic_subgradient_oracle} \& \ref{assume:channel_noise_independence}, and $\psi(\xv^*) \leq R^2$ for some $R$.
If
\ea{
\tau + 2\gamma = 1,
\label{eq:condition}
}
the expected suboptimality gap is bounded as
$$\esf(K) \leq \frac{20R\loge{K\sqrt{n}}}{K^{(1 - \gamma)/2}\sqrt{1 - \lambda}}\Paren{\Omega^2 + \frac{c_0^2\Delta^2d}{4} + \frac{c_0^2 \sigma^2 d}{2 \gamma c_1}}^{\f 12},$$
where $\lambda$ is the second largest eigenvalue of the matrix $\Pv$.
\end{theorem}

Note that $R$ is different from the rate $\Rsf$.
The suboptimality gap given by Thm. \ref{thm:main_convergence_result} is conditioned on the event that \gls{dlmd} succeeds.
The condition in \eqref{eq:condition} highlights the trade-off between the rate of decrease of confidence sequence and the increase in power allocation.
Note that a lower $\gamma$, i.e. faster convergence rate $\Ocal(K^{-(1-\gamma)/2})$ means a larger $\tau$, i.e. a larger transmit power requirement.
The next theorem states that with appropriate design choices, \gls{dlmd} can be made to succeed with a (high) probability $P_{\rm S}(K)$.\\

\vspace{-1mm}
\begin{theorem}
\label{thm:quantizer_unsaturation} [\textbf{Lower bound on quantizer unsaturation probability}]
Consider the same setting as in Thm. \ref{thm:main_convergence_result}.
For an execution of \gls{dlmd} till horizon $K$ using finite-range quantizers of resolution $\Delta$ and dynamic range $U$, with channel noise and stochastic subgradient oracle of variances $\sigma^2$ and $\Omega^2$ respectively, we have:
\begin{equation}
    P_{\rm S}(K) \geq \Paren{1 - \frac{1}{(U - U_{\min})^2}\Paren{\frac{c_0^2\Delta^2}{2} \hspace{-1mm} + \hspace{-1mm} \frac{c_0^2 \sigma^2}{2\gamma c_1} \hspace{-1mm} + \hspace{-1mm} \Omega^2}}^{2Kd \vert E \vert},
\end{equation}
where, $U_{\min} = L$ is the Lipschitz constant of each $f_i$, and $|E|$ denotes the number of edges in the connectivity graph.
\end{theorem}

Thm. \ref{thm:quantizer_unsaturation} shows that if we choose $U$, $\Rsf$ and $\Delta$ appropriately, $P_{\rm S}(K)$ can be made greater than $1 - \nu$ for any specified $\nu$.
The proofs of the above theorems, along with some further insights, are presented in Sec. \ref{sec:proofs_sketch}.

\subsection{System design implications}
\label{subsec:implications}
Note that for a given $\tau$, the power requirement for \gls{dlmd} at any node $i \in [n]$ can be expressed from \eqref{eq:power_requirement}, as:
\begin{equation}
    \label{eq:average_power_requirement}
    \Psf(K) = \frac{1}{K}\sum_{k \in [K]}\alpha(k)^2\hspace{-2mm}\sum_{j \in \Ncal(i)}\norm{\boldsymbol{\de}_{ij}(k)}^2 \leq \frac{ndc_1U^2}{\tau + 1}K^{\tau},
\end{equation}
where we use the inequality $\sum_{k \in [K]}k^{\tau} \leq K^{\tau+1}/\Paren{\tau + 1}$, the fact that the quantized differential satisfies $\norm{\boldsymbol{\de}_{ij}(k)}_{\infty} \leq U$, and $|\Ncal(i)| \leq n$ for any $i$.\footnote{$\Psf(K)$ is the power requirement at any particular node $i$.\IEEEauthorrefmark{4}}
In \eqref{eq:average_power_requirement}, the subscript $i$ is dropped for simplicity.\footnote{
In the more general case, if different nodes use different power control sequences $\{\alpha_i(k)\}$, we would have a \textit{per-node power requirement}.}

In any practical application scenario, the desired suboptimality gap $\epsilon$ and the failure probability threshold $\nu$ would be specified.
Given these specifications, one would be interested in determining  \gls{dlmd}  parameters necessary to meet them. 
In particular, we would need to determine the number of iterations $K$ so that $\esf(K) < \epsilon$ subject to $P_S(K) > 1 - \nu$.
The parameters that determine the performance of the algorithm are (i) the data rate $\Rsf$, (ii) the quantizer dynamic range $U$, and (iii) the confidence parameter $\gamma$.
Note that the other parameters in the expressions (except for the constants $c_0$ and $c_1$) are determined from the relations $\tau + 2\gamma = 1$, and $\Delta = 2U/\Paren{2^{\Rsf}-1}$.
The choice of these parameters is application specific.
A lower value of $\gamma$ would decrease the number of iterations required to reach the desired suboptimality gap $\epsilon$, but lower $\gamma$ implies larger $\tau$, which translates to a larger transmit power requirement for the algorithm.
Furthermore, for a given value of the data rate $\Rsf$, increasing the quantizer dynamic range $U$ not only implies a larger number of iterations to reach a given suboptimality gap $\epsilon$, but also a higher success probability $P_S(K)$.
For a given value of the dynamic range $U$, increasing $\Rsf$ (which translates to a higher network bandwidth requirement) and subsequently decreasing $\Delta$, implies higher $P_S(K)$ as well as a smaller number of iterations for reaching the desired suboptimality gap $\epsilon$.
The particular choice of parameters needs to be optimized for convergence speed, success probability and power-requirement.
It is also worth mentioning that in our simulations, we have not used any scheme to determine the optimal choice for the hyperparameters $\gamma$ and $\tau$.
The choice of $\gamma =0.5$ ensures least power requirement and the slowest convergence rate for the algorithm.
Similarly, $\gamma =0.1$ was chosen to provide a result close to the lower bound of the noiseless case.
$\tau$ and $\gamma$ are scalar values that satisfy a linear relation between them, i.e. $\tau + 2\gamma = 1$; only one of them can be varied independently.
Since we have a closed form expression for the upper bound on the expected suboptimality gap (Thm. 1), as well as an upper bound on the power requirement from Eq. (19), the (approximately) optimal value of $\gamma$ (and hence $\tau$), can be found in the range $\gamma \in (0, 0.5]$ by doing a binary search on this interval.

\section{Proofs Sketch}
\label{sec:proofs_sketch}

In this section, we present a sketch of the proofs of Lemma \ref{lem:convergence_rate_power_requirement_trade-off}, Thm. \ref{thm:main_convergence_result}, and Thm. \ref{thm:quantizer_unsaturation}.
These  three results characterize the performance of \gls{dlmd} in terms of conditions for convergence, suboptimality gap and probability of success, respectively.
A principle feature of the proof is that we ensure asymptotic consensus of the nodes' states despite the presence of two distinctive sources of noise as mentioned in Sec. \ref{subsec:our_contributions}.

Apart from the fact that the subgradient oracle has inherent stochasticity whose variance is assumed to be bounded according to Asm. \ref{assume:stochastic_subgradient_oracle}, stochasticity in the subgradient steps that \gls{dlmd} takes in line 17 of Alg.  \ref{algo:DLMD-DiffEx_pseudocode} also comes from quantization and channel noise.
The design of our finite dynamic range probabilistic quantizer ensures that the quantization noise has bounded variance when the conditions in Lemma \ref{lem:unbiased_bounded_variance_quantizer} are satisfied.
In order to bound the overall success probability, the major challenge is to ensure that the variance of the noise accumulated due to additive noisy channel remains bounded as the algorithm proceeds.
From \eqref{eq:noise_accumulation} (which is derived in Appendix \ref{appendix:proof_of_lemma_2}), we see that the noise variance grows linearly with the iteration index $k$.
If nothing is done to control this linear growth, the algorithm would diverge.

\gls{dlmd} makes use of \textit{power control} and \textit{diminishing confidence} sequences, respectively $\{\al(k)\}$ and $\{\beta(k)\}$, to ensure that the accumulated noise variance remains bounded as $k \to \infty$.
Lemma \ref{lem:convergence_rate_power_requirement_trade-off} gives a sufficient criteria which, if satisfied by these sequences, would ensure that the accumulated channel noise variance remains bounded.
The proof of Lemma \ref{lem:convergence_rate_power_requirement_trade-off} in Appendix \ref{appendix:proof_of_lemma_2} analyzes the discrepancy between the true state of a node i.e. $\zv_j(k)$ and the noisy version used by its neighbors for computing consensus, i.e. $\ytv_{ij}(k)$.
The discrepancy consists of the sum of two terms: (i) the quantization error, and (ii) the accumulated channel noise.
In other words, the effects of quantization noise and channel noise are simply additive and this observation simplifies the analysis a great deal.
From independence (Asm. \ref{assume:channel_noise_independence}), the variance of these two contributions can be added up and we focus on them separately.

The proof of Thm. \ref{thm:main_convergence_result} considers geometric sequences for $\alpha(k)$ and $\beta(k)$.
The requirement $\tau + 2\gamma = 1$ in the statement of the theorem is the sufficiency requirement of Lemma \ref{lem:convergence_rate_power_requirement_trade-off}. 
The proof first shows that the variance of the effective stochasticity in the subgradient step of \gls{dlmd}, i.e. $\xi^2$ is bounded.
The rest of the proof (presented in Appendix \ref{appendix:proof_of_thm_1}) is quite similar to that in \cite{duchi2012TAC}, and consists of the following steps:
\begin{enumerate}[1.]
    \item Deriving an upper bound for the network consensus error $\norm{\zvo(k) - \zv_i(k)}_*$ in the presence of confidence sequence $\beta(k)$, which we do in Lemma \ref{lem:convergence_doubly_stochastic_matrices} and Lemma \ref{lem:consensus_error_upper_bound}.
    \item Obtaining an upper bound on the expected suboptimality gap $\Expect f(\xvh_i(K)) - f(\xv^*)$ in terms of $\Expect \norm{\zvo(k) - \zv_i(k)}_*$.
    \item Substituting the expression for $\norm{\zvo(k) - \zv_i(k)}_*$ and getting the result.
\end{enumerate}

Finally, Thm. \ref{thm:quantizer_unsaturation} gives an expression for the \textit{success probability} of \gls{dlmd}, i.e. that none of the quantizers get saturated.
The proof of this (presented in Appendix \ref{appendix:proof_of_theorem_2}) computes the $\ell_{\infty}$ norm of the differential vector being exchanged at any iteration, and makes use of Chebyshev's inequality to show that for a high enough dynamic range and data rate, it is relatively unlikely that the quantizers would ever saturate.

\section{Numerical simulations}
\label{sec:numerical_simulations}

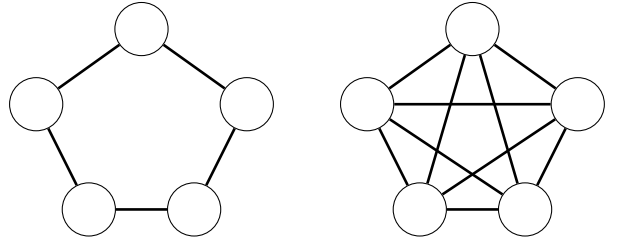
\begin{figure}[h!]
    \centering
    \begin{tikzpicture}[scale = 0.4]
  
        \node (one) [draw, circle,inner sep=0.25cm] at (11,0){};
        \node (two) [draw, circle,inner sep=0.25cm] at (14.5,2.5){};
        \node (three) [draw, circle,inner sep=0.25cm] at (18,0){};
        \node (four) [draw, circle,inner sep=0.25cm] at (16.25,-3.5){};
        \node (five) [draw, circle,inner sep=0.25cm] at (12.75,-3.5){};
        
        \draw[line width=1pt] (one) -- (two);
        \draw[line width=1pt] (one) -- (three);
        \draw[line width=1pt] (one) -- (four);
        \draw[line width=1pt] (one) -- (five);
        \draw[line width=1pt] (two) -- (three);
        \draw[line width=1pt] (two) -- (four);
        \draw[line width=1pt] (two) -- (five);
        \draw[line width=1pt] (three) -- (four);
        \draw[line width=1pt] (three) -- (five);
        \draw[line width=1pt] (four) -- (five);
  
        \node (one) [draw, circle,inner sep=0.25cm] at (0,0){};
        \node (two) [draw, circle,inner sep=0.25cm] at (3.5,2.5){};
        \node (three) [draw, circle,inner sep=0.25cm] at (7,0){};
        \node (four) [draw, circle,inner sep=0.25cm] at (5.25,-3.5){};
        \node (five) [draw, circle,inner sep=0.25cm] at (1.75,-3.5){};
        
        \draw[line width=1pt] (one) -- (two);
        \draw[line width=1pt] (two) -- (three);
        \draw[line width=1pt] (three) -- (four);
        \draw[line width=1pt] (four) -- (five);
        \draw[line width=1pt] (five) -- (one);
        
    \end{tikzpicture}
    \caption{Ring with $2$ neighbors (left) and fully-connected (right) topologies over a network of $n=5$ nodes.}
    \label{fig:network_topologies}
\end{figure}

In this section, we evaluate the performance of \gls{dlmd} for a common learning scenario: Training a {Support Vector Machine (SVM)} over data distributed across multiple nodes of a network.
As depicted in  Fig. \ref{fig:network_topologies}, we consider two topologies: A 2-neighbour ring network and a fully-connected network. 
In both cases, we choose  $n = 10$ nodes.  

The local datasets for SVM training (as discussed in Remark \ref{rem:stochastic optimization}) are as follows: 
Node $i \in [n]$ has access to the data points $\{\av_{ij}, b_{ij}\}$, for $ j \in [m]$, where $m$ is the local dataset size, $\av_{ij} \in \Real^d$ with $d = 30$, and $b_{ij} \in \{-1, +1\}$ is the label identifying the class that datapoint $\av_{ij}$ belongs to. 
We consider $m = 10$ datapoints per node. Each point is generated from Gaussian distributions with different means for each of the two classes: $b_{ij}=-1$ and $b_{ij}=+1$.

In addition, we also consider the case in which the distribution of datapoints is polarized, that is a particular node contains datapoints from one particular class only.
The global objective is to train an SVM which can correctly classify the entire dataset (across all the nodes).
Such a polarized data distribution prevents local training from being successful and the SVM must be trained in a decentralized fashion.
Decentralized training of a soft-margin SVM can be formulated according to \eqref{eq:distributed_optimization_problem_statement} as the following $\ell_2^2$-regularized hinge loss minimization problem:
\begin{equation}
    \label{eq:SVM_optimization}
    \hspace{-2mm}\minimize_{\xv \in \Real^d} \frac{1}{n}\hspace{-1mm}\sum_{i \in [n] }\Br{\frac{1}{m}\hspace{-1mm}\sum_{j \in [m]}\hspace{-1mm}\max\lcb {0,1-b_{ij}\xv^\top\av_{ij}} \rcb + \frac{\mu}{2}\norm{\xv}_2^2},
\end{equation}
where, $\mu = 0.1$ is the regularization parameter.
In other words, the solution of \eqref{eq:SVM_optimization} yields the linear separator which can classify these $N = nm$ datapoints with the maximum margin. 

For this optimization problem, we wish to characterize the performance of \gls{dlmd} in two ways: 
$(i)$ First, we plot the suboptimality gap vs. number of iterations for different choices of confidence and power control parameters, $\gamma$ and $\tau$.
From Thm. \ref{thm:main_convergence_result}, lower values of $\gamma$ imply faster convergence, but the trade-off \eqref{eq:condition} requires us to have a larger $\tau$, i.e. a larger power requirement.
$(ii)$ Next, we investigate the suboptimality gap vs. success probability trade-off and comment on how the choice of the dynamic range $U$ affects this trade-off. 
A summary of the parameters used for simulations is presented in Table \ref{tab:simulation_parameters}.
The proximal function is taken to be the quadratic function, $\psi(\xv) = \frac{1}{2}\norm{\xv}_2^2$.
Since the optimal step-size depends on unknown parameters, we choose any small step size which gives a reasonable convergence.
However, once chosen, it is kept unchanged throughout the different simulations to maintain consistency of comparison.

For fully-connected networks, the adjacency matrix $\Pv \in \Real^{n \times n}$ is given by $\Pv_{ij} = 1/n$ for all $i,j \in [n]$.
For a ring network with $2$-neighbors, $\Pv$ is given by:
\begin{equation}
    \Pv_{ij} = 
    \begin{cases}
      1/3 & j = i \pm 1 \mod n\\
      0 & \text{otherwise},
    \end{cases}
\end{equation}

Fully-connected networks have a larger spectral gap (i.e. smaller $\lambda$) than ring-networks and hence achieve consensus faster (which depends inversely on the spectral gap).

\begin{table}
  \centering
   \caption{Parameters for simulations in Sec. \ref{sec:Duality gap vs. number of iterations} and \ref{sec:Duality gap vs. success probability}.}
  \begin{tabular}{|P{1.5cm}|P{1.5cm}|P{1.5cm}|}
    \hline
    Figures & Parameter  & Value \\ \hline
    \multirow{3}{*}{All Figs.}     & n                              & 10    \\
                         & d                              & 30    \\
                         & m                              & 10    \\ \hline
    \multirow{3}{*}{Figs. \eqref{fig:suboptimality_gap_no_of_iterations}} & $\Rsf$                              & 6     \\
                         & $\sgs$                         &  0.1     \\
                         & U                              &   100    \\ \hline
    \multirow{4}{*}{Figs. \eqref{fig:success_dynamic_range_variation}} & $\Rsf$                              & 3     \\
                         & $\sgs$                         &  0.05    \\
                         & U                              &   0.8 \text{to} 1.8    \\
                         & K                              &   75 \\ \hline
  \end{tabular}
  \label{tab:simulation_parameters}
\end{table}

\subsection{Suboptimality gap vs. number of iterations}
\label{sec:Duality gap vs. number of iterations}
\begin{figure*}[h!]
 \begin{subfigure}[b]{0.5\textwidth}
         \centering
         \includegraphics[width=1.1\textwidth]{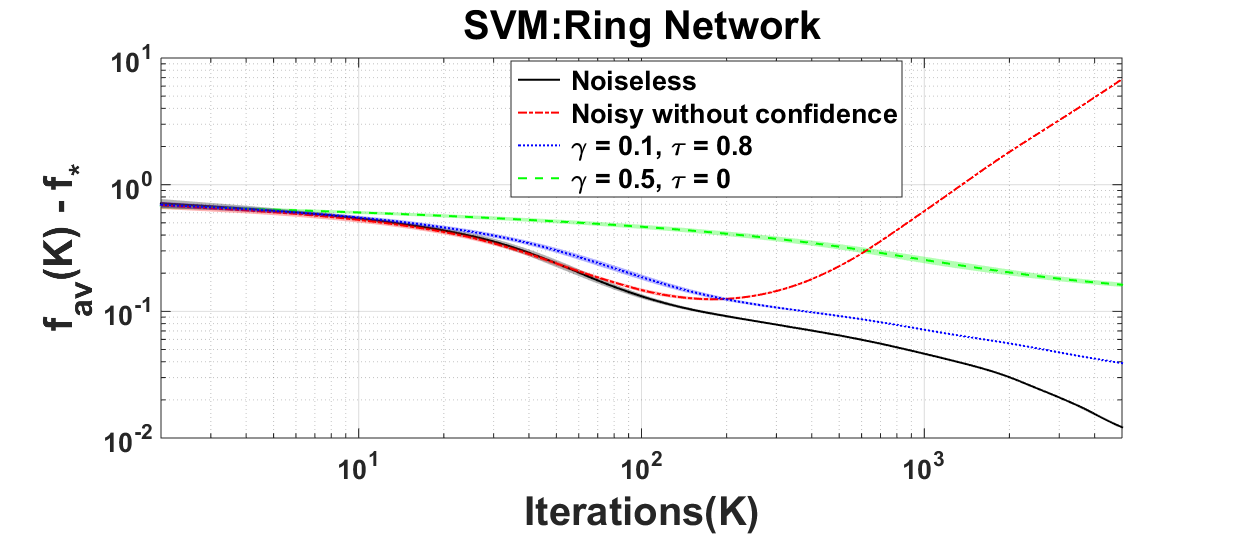}
     \end{subfigure}
     \hfill
     \begin{subfigure}[b]{0.5\textwidth}
         \centering
         \includegraphics[width=1.1\textwidth]{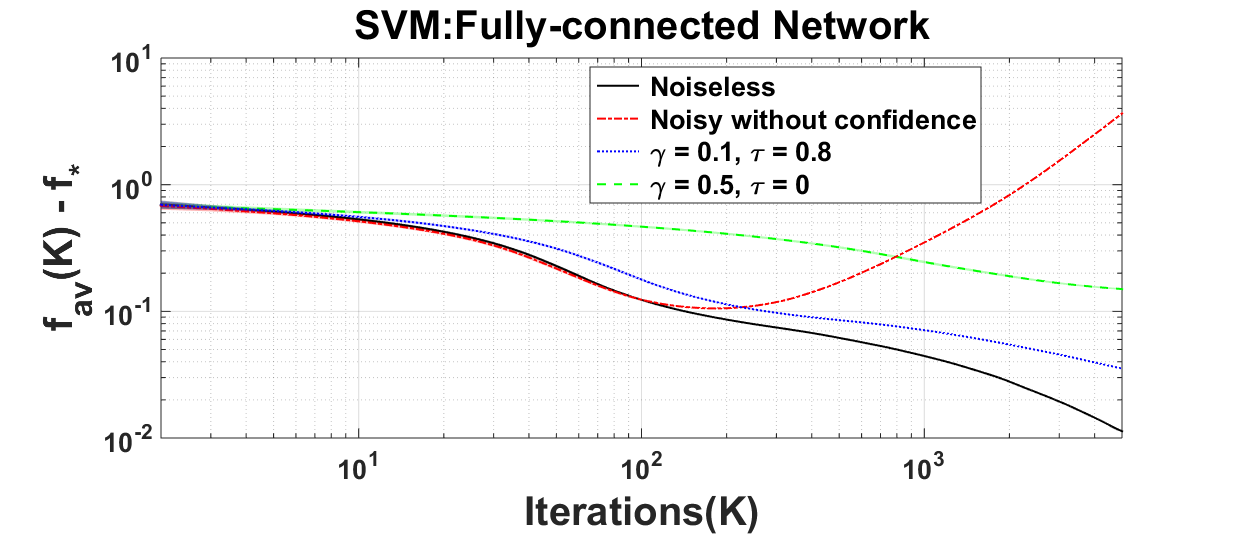}
     \end{subfigure}
          \caption{Suboptimality gap vs. number of iterations in Sec. \ref{sec:Duality gap vs. number of iterations} for the two topologies in Fig. \ref{fig:network_topologies}.
          }
     \label{fig:suboptimality_gap_no_of_iterations}
\end{figure*}

In Fig. \ref{fig:suboptimality_gap_no_of_iterations}, we plot the suboptimality gap across nodes, i.e.,
\ea{
f_{av}(K) - f_* =  \f 1 n \sum_{i\in [n]} f(\xv_i(k)) -f(\xv^*),
\label{eq:sub 2}
}
versus the number of iterations $K$.\footnote{Note that \eqref{eq:sub 2} is the empirically averaged suboptimality gap across nodes for a given $K$, while \eqref{eq:suboptimality_gap} is the expected gap of one node across iterations over sources of stochasticity like channel noise, quantization and subgradient oracle.}
Here, $f_* = f(\xv^*)$.

Fig. \ref{fig:suboptimality_gap_no_of_iterations} has four plots.
(i) The label \textit{``noiseless''} (black) refers to the setting in which the noise variance is set to zero, i.e. $\sgs=0$. Accordingly, the confidence and  power control sequences are set to the all-one sequence. 
In this case, \gls{dlmd} reduces to the naive DLMD algorithm of Sec. \ref{subsec:Decentralized_Lazy_Mirror_Descent} but with an additional noise term arising from the finite-range quantization. 
Given the noiseless setting, this curve provides a lower bound for the other simulations in the plot.
(ii) The label \textit{``noisy without confidence''} (red) considers the noisy setting but we use naive DLMD as in \cite{duchi2012TAC} without any \textit{confidence} or \textit{power control}, i.e. once again, they are set to all-ones sequence.
In this case, we see that the algorithm diverges when nothing is done to prevent channel noise accumulation.
The blue and green plots correspond to \gls{dlmd} with two different choices of parameters.
The channel noise variance in these cases is the same as the \textit{noisy} plot; yet \gls{dlmd} achieves convergence under this setting.
$\gamma = 0.5$ shows slower convergence than $\gamma = 0.1$, which is consistent with Thm. \ref{thm:main_convergence_result}.
Note that smaller $\gamma$ implies larger $\tau$ (according to \eqref{eq:condition}) and a subsequently larger power requirement for \gls{dlmd}.
Each plot is obtained by averaging over $5$ different realizations.
The initialization is kept the same in all cases for a fair comparison. 
The data rate and the dynamic range for this set of simulations is kept high enough so that the quantizers never saturate.

\subsection{Suboptimality gap vs. success probability}
\label{sec:Duality gap vs. success probability}

\begin{figure*}[h!]
\captionsetup[subfigure]{justification=centering}
    \begin{subfigure}[b]{0.5\textwidth}
         \centering
         \includegraphics[width=1.1\textwidth]{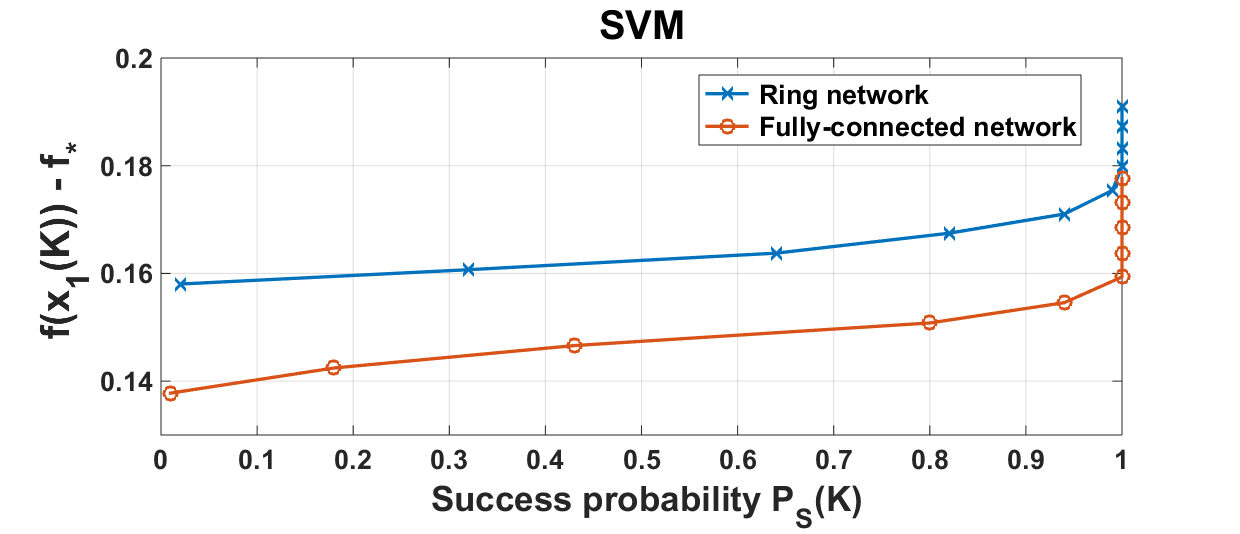}
         \caption{Node $1$ suboptimality vs. success probability.}
    \label{fig:success_suboptimality_curve}
     \end{subfigure}
     \hfill
     \captionsetup[subfigure]{justification=centering}
     \begin{subfigure}[b]{0.5\textwidth}
         \centering
         \includegraphics[width=1.1\textwidth]{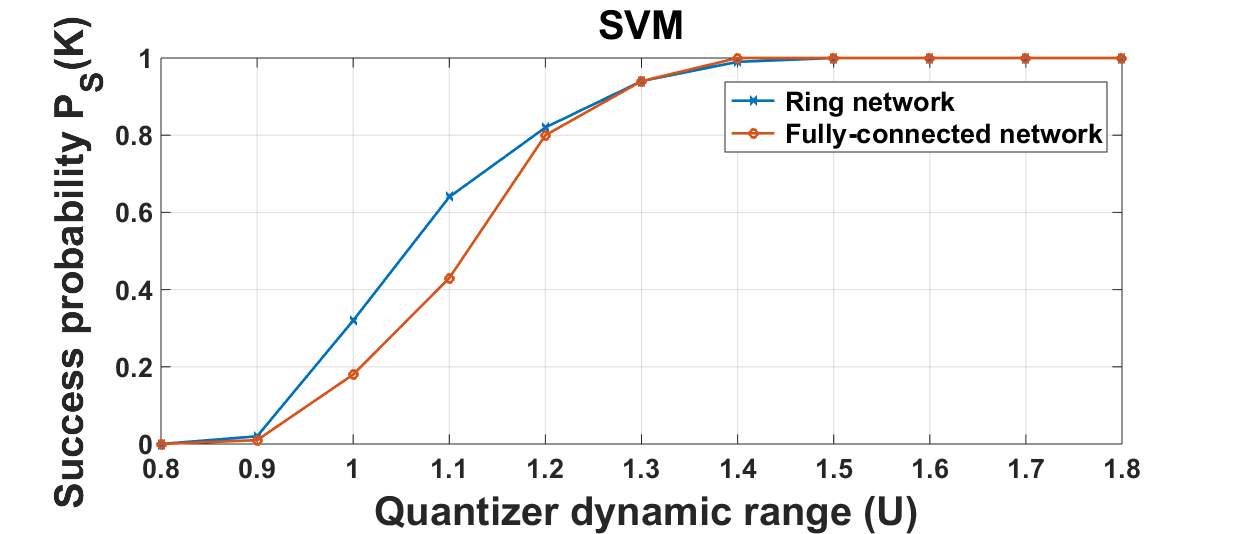}
        \caption{Variation of success probability with quantizer dynamic range.}
    \label{fig:dynamic_range_variation}
    
     \end{subfigure}
          \caption{Suboptimality gap vs. probability of success and Success Probability vs. Dynamic Range in Sec. \ref{sec:Duality gap vs. success probability} for the two topologies in Fig. \ref{fig:network_topologies}.}
   \label{fig:success_dynamic_range_variation}
\end{figure*}

\begin{figure*}[h!]
\captionsetup[subfigure]{justification=centering}
 \begin{subfigure}[b]{0.5\textwidth}
         \centering
         \includegraphics[width=1.1\textwidth]{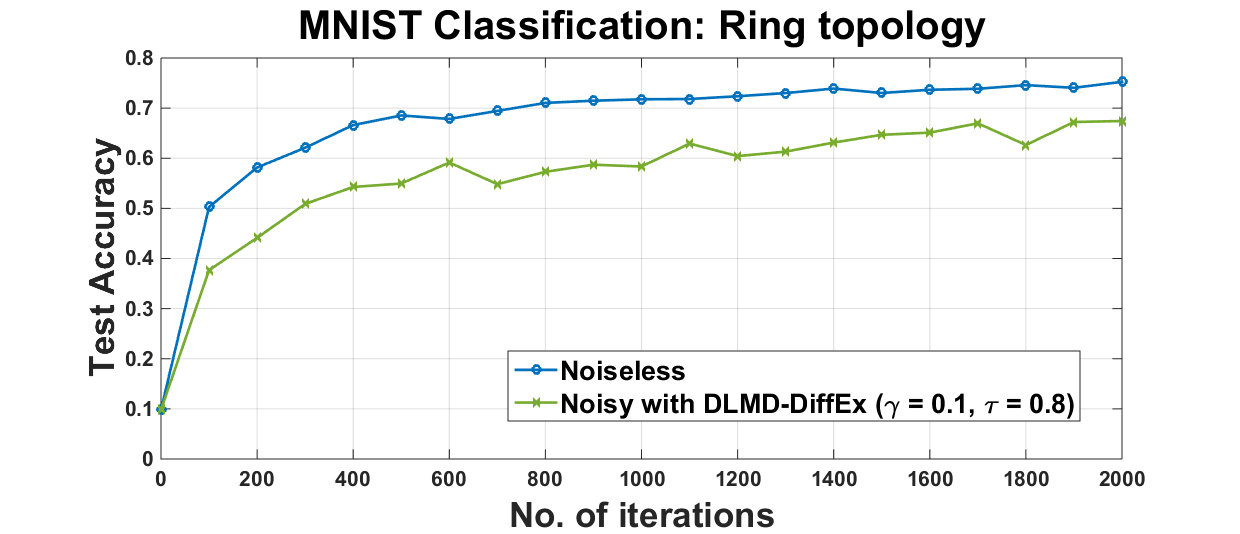}
         \caption{Test accuracy using fully-connected neural network}
    \label{fig:test_acc_NN}
     \end{subfigure}
     \hfill
     \captionsetup[subfigure]{justification=centering}
     \begin{subfigure}[b]{0.5\textwidth}
         \centering
         \includegraphics[width=1.1\textwidth]{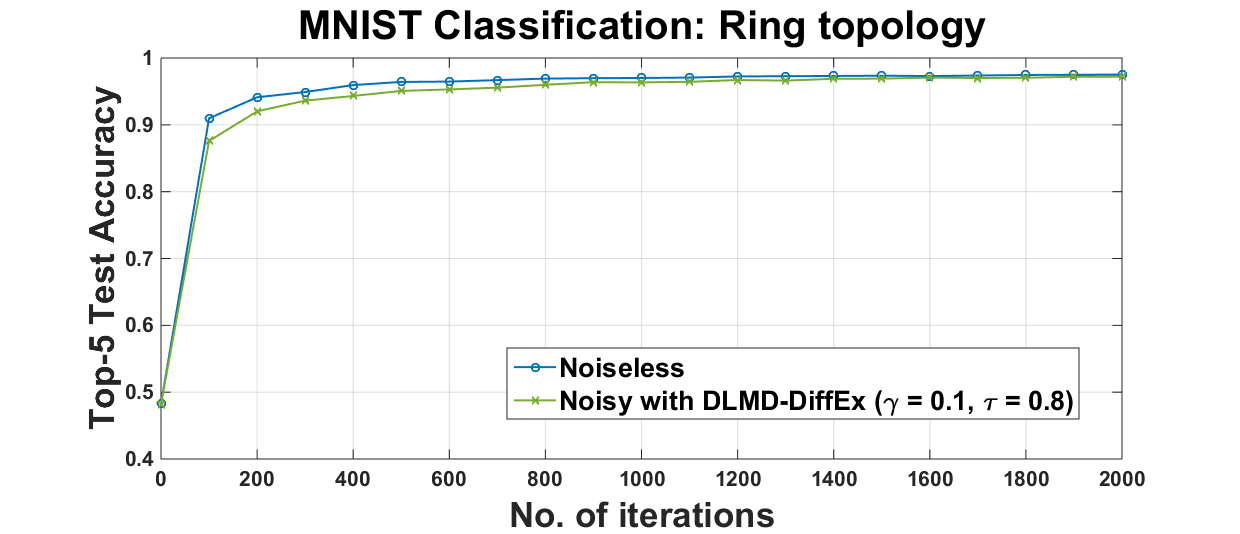}
        \caption{Top-$5$ test accuracy using fully-connected neural network}
    \label{fig:test_kbest_NN}
    
     \end{subfigure}
          \caption{Performance of a \textsc{DLMD-DiffEx} for decentralized training of a fully-connected neural network}
\end{figure*}

In Fig. \ref{fig:success_suboptimality_curve}, we empirically demonstrate the trade-off between the success probability $P_S(K)$ and the suboptimality gap at a particular node, that is $f(\xv_1(k)) -f(\xv^*)$.
%
%
The plots correspond to two different network topologies, and for the same value of $P_S(K)$, they have different suboptimality gaps because of the difference in their spectral gaps.

We next turn our attention to Fig. \ref{fig:dynamic_range_variation} which plots the variation of the success probability with the quantizer dynamic range $U$.
From this figure, we notice that as $U$ increases, $P_S(K)$ also increases according to Thm. \ref{thm:quantizer_unsaturation}.
This implies that $\Delta = 2U/\Paren{2^{\Rsf} - 1}$ increases, subsequently increasing the expected suboptimality gap according to Thm. \ref{thm:main_convergence_result}.
For a specific value of $U$, \gls{dlmd} is executed $100$ times and we count the number of times the algorithm execution is completed up to $K= 75$ iterations (reaching which without the quantizers getting saturated is deemed as a success).
The final suboptimality gap of Fig. \ref{fig:success_suboptimality_curve} is obtained by averaging out only over those successful iterations.
Note that after a certain value of $U$, $P_S(K)$ reaches $1$.
This occurs because there is a value of $U$ after which success is practically guaranteed and increasing $U$ beyond that point only deteriorates the suboptimality gap.
To move along either of the plots in Figs. \ref{fig:success_suboptimality_curve} and \ref{fig:dynamic_range_variation}, the dynamic range $U$ of the finite-rate quantizers is varied.

Another worthwhile thing to note in Fig. \ref{fig:dynamic_range_variation} is the fact that for the same value of $U$, the ring-network has higher success probability than the fully-connected network.
This is consistent with the fact that $|E|$, the number of edges, appears as an exponent in the expression for $P_S(K)$ in Thm. \ref{thm:quantizer_unsaturation}.
Moreover, $\lambda$ does not appear in the expression for success probability in Thm. \ref{thm:quantizer_unsaturation}.

The simulation outputs are consistent with the system design implications in Sec. \ref{subsec:implications}.
\gls{dlmd} requires us to look into transmit power consumption and data rate requirements for any specified application.
Note that if for a given choice of parameters, \gls{dlmd} has a failure probability $1 - P_S(K)$, repeating the algorithm execution $M$ times reduces the failure probability to $(1 - P_S(K))^M$, which can be made arbitrarily small by making $M$ appropriately large.
However, repeating the algorithm execution multiple times increases the wall-clock time elapsed for finding the optimal solution by a factor of $M$ and hence, it is important to take into account how large $P_S(K)$ should be while choosing algorithm parameters $\gamma$ and $U$.
For time-critical applications, where it undesirable to repeat the algorithm execution several times, one would want $P_S(K)$ to be relatively large; in other words, it is desirable to operate at a point on the far right of Fig. \ref{fig:success_suboptimality_curve}.

\subsection{MNIST digit Classification using Fully-Connected Neural Networks (NN)}
\label{sec:simulations_non_convex}

In this section, we consider the problem of classifying handwritten digits from the MNIST dataset \cite{lecun_mnist}.
We consider a ring network comprised of $5$ nodes.
The data is distributed in way so that each node consists of data from only two of the ten possible classes.
This once again means that if each of the nodes were to independently train a model without collaborating with others, they would not have been successful.
We consider a fully-connected network for the classification task.
The network takes in a $28 \times 28$ dimensional image and flattens it into a vector.
This is followed by a dense layer with $64$ nodes, finally followed by another dense layer with $10$ nodes.
The proximal function is taken to be $\frac{1}{2}\norm{\xv - \xv_{init}}_2^2$, where $\xv_{init}$ is the random initialization.

In Figs. \ref{fig:test_acc_NN} and \ref{fig:test_kbest_NN}, we plot the test accuracy and the top-$5$ accuracy of the neural network trained in a decentralized fashion.
Top-$5$ accuracy refers to if the ground truth of a test input belongs to the top $5$ predicted classes.
We consider a dynamic range of $U = 30$ for our finite range quantizers, and $M = 2^{\Rsf} = 100$ quantization points per dimension.
The noise variance is taken to be $\sigma^2 = 0.1$.

Although we have done the analysis only for convex objective functions, we see that our algorithm converges for non-convex objectives (i.e. NN training) as well under rate-constrained and noisy settings.
The blue plots correspond to the accuracy of the noiseless setting.
The green plots are obtained when the differential exchanges are noisy, and \textsc{DLMD-DiffEx} is used as the optimization algorithm with a confidence parameter of $\gamma = 0.1$, and a power control parameter of $\tau = 0.8$.
We leave the analysis for non-convex settings for future research.

\vspace{-4mm}
\section{Additional Remarks}
\label{sec:additional_remarks}
Before concluding the paper, we mention some limitations of our problem formulation and possible future research directions.
We note that our channel model incorporating data-rate constraints and additive channel noise, which gives rise to the associated trade-offs as uncovered in this paper, can also be extended to other classes of decentralized optimization algorithms such as primal-only and second order methods.
Moreover, the corresponding trade-offs in the study of non-convex settings, for example decentralized training of neural networks, is still open.

In particular, using Lazy mirror descent for non-convex objectives requires the use of an appropriate proximal function $\psi(\cdot)$ that makes effective use of the geometry of the objective function and with respect to which projection can be computed.
For our simulations in Sec. \ref{sec:simulations_non_convex}, we have used a quadratic proximal function which reduces lazy mirror descent to the standard subgradient algorithm.
For primal-only algorithms, works like \cite{gorbunov2021marina} study communication compression for non-convex optimization.
The extension of their proposed algorithm MARINA to take into account additive noise accumulation is worth considering. 

\textbf{Channel fading}: A practical concern for wireless environments is fading.
We would like to remark that in this work, we have considered perfect channel knowledge at the receiving node, so that the effects of fading can be inverted.
This requires us to deviate from our choice of confidence and power control sequences as geometric sequences, in order to effectively incorporate channel inversion techniques.
It is also possible to extend the framework proposed in Sec. \ref{sec:system_model} to  the case of heterogeneous links by suitably choosing the confidence and power sequences to be dependent on the nodes.

\textbf{Random link failures}: Another relevant concern for wireless networks is that the links may be randomly time-varying.
\cite{duchi2012TAC} considers the case of stochastic communication where the only assumption made on the randomly varying network topology is that the expected spectral gap is constant, i.e. $\lambda_2(G) = \lambda_2(\mathbb{E}[\Pv(t)^{\top}\Pv(t)])$ is independent of $t$.
It is possible to extend the proof provided in \cite[Sec. VII]{duchi2012TAC} which considers a sequence $\{\Pv(t)\}_{t=0}^{\infty}$ of random matrices that describe the randomly varying network topology typical of wireless networks, and utilize properties of Markov chain mixing for completing the proof.

\textbf{Extension to compression schemes}: A possible extension of the proposed algorithm is by incorporating more general compression schemes \cite{Koloskova2020Decentralized}, by deriving a result similar to Lemma 1 for general compression operators, as such operators have bounded error.

For general compression operators such as the ones considered in 
\cite{beznosikov2020biased, qian_prox_SGD_RDA} (which include biased compressors), the variance of the quantization error is correlated with the value of the input to the quantizer.
These works consider compressors $\Ccal(\cdot)$ characterized by a parameter $0 < \delta \leq 1$ such that $\Expect\norm{\xv - \Ccal(\xv)}^2 \leq (1- \delta)\norm{\xv}^2$.
To ensure that our analysis is tight, the fact that this upper bound on the variance depends on the actual value of the parameter might necessitate a slightly different analysis.
This is because in our work, the expected quantization error is upper bounded by a term that only requires knowledge of the quantizer dynamic range $U$. 

Lastly, it is interesting to note that the transmit power is dependent on $K$ \eqref{eq:average_power_requirement}.
The setting in which the transmit power is constant is obtained by setting $\tau = 0$.
%
This will require $\gamma = 1/2$ (according to eq. (17)), which means that the convergence rate drops to $\Ocal(K^{-1/4})$.
If we consider an iteration horizon $K$ beforehand for which we decide we will execute the algorithm, eq. (19) tells us the power requirement of the algorithm.
It is still an open question whether this tradeoff between the power requirement and the convergence rate is Pareto optimal.
An intuitive way to understand this dependence of power on $K$ is to note that for non-smooth functions, it is possible that the subgradient at each step has equal magnitude, i.e. it is possible that the update $\zv(k+1) = \zv(k) + \gv(k)$ in the subgradient step of lazy mirror descent is equally informative for all $k$.
Note that since noise accumulates, for $\gv(k)$ to have finite variance for large values of $k$,  it becomes important to geometrically increase power in the later stages of iteration, and hence the dependence of the transmission power on $K$ appears.
However, for smooth and strongly convex objectives, the magnitude of the gradients decrease as iterations $K$ increase.
For such objective functions, it might be possible that the power requirement is independent of $K$.
The case for strongly convex and smooth objectives is yet to be studied.

\section{Conclusions}
\label{sec:conclusions}

In this work, we considered the problem of decentralized optimization over a network of nodes, where the information exchange between them is constrained by network connectivity and finite data rate, as well as corrupted by an additive channel noise.
We argued that the simultaneous presence of both finite data rate and channel noise gives rise to a set of difficulties that do not arise when either of them is considered individually.
In particular, we show that straightforward extensions of existing algorithms to such settings can lead to noise accumulation and subsequent divergence.
For this reason, we propose a novel algorithm: \textit{Decentralized Lazy Mirror Descent with Differential Exchanges} (\gls{dlmd}).
We showed that the algorithm parameters in \gls{dlmd} allow us to trade off power consumption, data rate requirements, convergence rate and success probability.
We analyzed \gls{dlmd} and demonstrated that our numerical simulations were consistent with our theoretical analysis.

\bibliographystyle{IEEEtran}
\bibliography{DLMD_DiffEx}

\appendices

\section{Proof of Lemma \ref{lem:convergence_rate_power_requirement_trade-off}}
\label{appendix:proof_of_lemma_2}

From line 8 of Alg. \ref{algo:DLMD-DiffEx_pseudocode}, note that the quantized state differential is obtained as
\ea{
\boldsymbol{\delta}_{ij}(k) = \zv_j(k) - \yv_{ij}(k-1) + \Deltav_{ij}(k),
}
where $\Deltav_{ij}(k)$ is the quantization error.
Line 13 in Alg. \ref{algo:DLMD-DiffEx_pseudocode} gives 
\ea{
\deltatv_{ij}(k) = \alpha(k)^{-1}\rv_{ij}(k) = \deltv_{ij}(k) + \Tilde{\nv}_{ij}(k),
}
where $\Var\Paren{\Br{\Tilde{\nv}_{ij}(k)}_s} = \sigma^2/\alpha(k)^2$ for $s \in [d]$.
The running estimate of node $j$'s state maintained by node $i$ is updated as:
\begin{align}
    \label{eq:running_estimate}
    &\ytv_{ij}(k) = \tilde{\deltv}_{ij}(k) + \ytv_{ij}(k-1) \\
    &= \zv_j(k) + \Deltav_{ij}(k) + \ntv_{ij}(k) + \Br{\ytv_{ij}(k-1) - \yv_{ij}(k-1)}. \nonumber
\end{align}

To track  the evolution of the discrepancy $\ytv_{ij}(k-1) - \yv_{ij}(k-1)$ in $k$, note that the proxy at node $j$ is updated as $\yv_{ij}(k) = \deltv_{ij}(k) + \yv_{ij}(k-1)$, and the estimate at node $i$ is:
\begin{align*}
    \ytv_{ij}(k) & = \deltatv_{ij}(k) + \ytv_{ij}(k-1) \\ 
&= \deltv_{ij}(k) + \tilde{\nv}_{ij}(k) + \ytv_{ij}(k-1). \nonumber
\end{align*}
Accordingly the discrepancy above grows  as
\begin{equation}
    \label{eq:discrepancy_evolution_recursion}
    \ytv_{ij}(k) - \yv_{ij}(k) = \ytv_{ij}(k-1) - \yv_{ij}(k-1) + \ntv_{ij}(k)
\end{equation}

Unrolling this telescopic sum gives $\ytv_{ij}(k) - \yv_{ij}(k) = \sum_{l\in [k]}\ntv_{ij}(l)$.
From Asm. \ref{assume:channel_noise_independence}, the variance of the accumulated noise (which adds to the stochasticity of the subgradient as seen in \eqref{eq:dlmd_update_equation_effective_stochasticity}) is
\begin{align}
\label{eq:accumulated_noise_variance_general}
     &\Var\Paren{\Br{\ntv_i(k)}_s} =  \Var\Paren{\Br{\sum_{j  \in [n] \setminus \{i\} } \sum_{l\in [k]}\Wv_{ij}(k)\ntv_{ij}(l)}_s} \nonumber \\
     &\leq \sum_{j  \in [n] \setminus \{i\} }\sum_{l\in [k]}\Var\Paren{\Br{\beta(k)\Pv_{ij}(k)\ntv_{ij}(l)}_s} \nonumber \\
     &\leq \sum_{j  \in [n] \setminus \{i\} }\sum_{l\in [k]}\frac{\beta(k)^2}{\al(l)^2}\Pv_{ij}^2 \sigma^2 \nonumber \\ 
     & \leq \sgs \beta(k)^2 \sum_{l\in [k]} \frac{1}{\al(l)^2}.
\end{align}
From \eqref{eq:accumulated_noise_variance_general}, we see that the condition $\sum_{l\in [k]}\al(l)^{-2} = \Theta\Paren{\beta(k)^{-2}}$ in \eqref{eq:condition al be} is sufficient to guarantee that this variance stays bounded as $k \to \infty$. 
Note that taking a summation over all nodes weighted by $\Pv_{ij}$ is equivalent to taking consensus over the neighborhood. 

\section{Proof of Theorem \ref{thm:main_convergence_result}}
\label{appendix:proof_of_thm_1}

The proof of Thm. \ref{thm:main_convergence_result} closely follows the convergence proof of distributed dual averaging \cite{duchi2012TAC} with appropriate modifications to take into account the effects of additional sources of stochasticity, i.e. quantization and channel noise.
The proof structure is as follows:
(i) We begin by introducing some necessary lemmas from \cite{duchi2012TAC} which we state without proof.
(ii) Then,  we show that as a consequence of Lemma \ref{lem:unbiased_bounded_variance_quantizer}, Asm. \ref{assume:stochastic_subgradient_oracle}, Asm. \ref{assume:channel_noise_independence}, and with the {confidence} and {power control} sequences chosen to be geometric sequences, the effective stochasticity in the subgradient vector due to stochastic oracle, quantization and channel noise has bounded variance.
(iii) Next, in Appendix \ref{subsec:upper_bound_consensus_error}, we derive an upper bound to the consensus error $\norm{\zvo(k) - \zv_i(k)}_*$ in the presence of \textit{diminishing confidence} but without additive channel noise or rate constraint. 
(iv)  Finally, in Appendix \ref{subsec:expected suboptimality_gap_general}, we combine the result in (iii) with that of (ii) to bound the overall suboptimality gap and conclude the proof.

\begin{lemma}{\bf \cite[Lemma 2]{duchi2012TAC}}
\label{lem:projection_lipschitz_continuity}
    For any arbitrary pair $\uv,\vv \in \Real^d$, it holds that  
    \ea{
    \norm{\Pi_{\Xcal}^{\psi}\Paren{\uv, \eta} - \Pi_{\Xcal}^{\psi}\Paren{\vv, \eta}} \leq \eta\norm{\uv - \vv}_*.}
\end{lemma}
\begin{lemma}{\bf \cite[Lemma 3]{duchi2012TAC}}
\label{lem:projection_sequence_lemma}
    Let $\left\{ \gv(k) \right\}_{k  \in \Nbb} \in \Real^d$ be an arbitrary sequence of vectors and consider the sequence $\{\xv(k)\}_{k\in \Nbb}$ defined by $\xv\Paren{k+1} = \mathlarger{\Pi}_{\xv \in \Xcal}^{\psi}\Paren{\sum_{s\in [k]}\gv\Paren{s}, \eta(k)}$.
    Then for any non-increasing sequence $\{\eta(k)\}_{k  \in \Nbb}$ of positive step-sizes, and for any $\xv^* \in \Xcal$, we have
    \ea{
        & \sum_{k \in [K]}\inprod{\gv\Paren{k}, \xv\Paren{k} - \xv^*}\leq \\  
        & \quad \quad \f 12\sum_{k \in [K]}\eta\Paren{k-1}\norm{\gv(k)}_*^2  +\frac{1}{\eta(K)}\psi\Paren{\xv^*}.
        \nonumber
    }
\end{lemma}
\begin{lemma}{\bf \cite[Lemma 4]{duchi2012TAC}}
\label{lem:centralized_individual_sequence_relation}
For the sequence $\{\yv(k)\}_{k  \in \Nbb}$ obtained as $\yv\Paren{k} = \Pi_{\xv \in \Xcal}^{\psi}\Paren{\zvo(k), \eta(k)}$, where $\zvo(k) = \frac{1}{n}\sum_{i \in [n] }\zv_i(k)$ is the node-averaged state sequence, we have
\begin{align}
    & \sum_{k \in [K]}\Br{f(\xv_i(k)) - f(\xv^*)} \leq 
    \label{eq:centralized_individual_sequence_relation} \\ 
    & \quad \quad \sum_{k \in [K]}\Br{f(\yv(k)) - f(\xv^*)} + L\sum_{k \in [K]}\hspace{-1mm}\eta(k)\norm{\zvo(k) - \zv_i(k)}. \nonumber
\end{align}
Furthermore, if the iteration-averaged iterates are denoted as $\yvh(K) = K^{-1}\sum_{k \in [K]}\yv(k)$ and, $\xvh_i(K) = K^{-1}\sum_{k \in [K]}\xv_i(k)$, we also have $ f(\xvh_i(K)) - f(\xv^*) \leq f(\yvh(K)) - f(\xv^*) + LK^{-1}\sum_{k \in [K]}\eta(k)\norm{\zvo(k) -\zv_i(k)}$.
\end{lemma}

From \eqref{eq:running_estimate} and \eqref{eq:discrepancy_evolution_recursion}, the running estimate of node $j$'s state at node $i$ is given as $\ytv_{ij}(k) = \zv_j(k) + \Deltav_{ij}(k) + \sum_{l\in [k]}\ntv_{ij}(l)$.
At any iteration $k$, the state update at node $i$ is then obtained as
\begin{align}
    \label{eq:dlmd_update_equation_effective_stochasticity}
    &\zv_i\Paren{k+1} = \Wv_{ii}(k)\zv_i(k) + \nonumber \\  &\sum_{j  \in [n] \setminus \{i\}}\hspace{-2mm}\Wv_{ij}(k)\Paren{\zv_j(k) + \Deltav_{ij}(k) + \sum_{l\in [k]}\ntv_{ij}(l)} + \gtv_i(k) \nonumber \\
    & = \sum_{j \in [n] }\Wv_{ij}(k)\zv_j(k) + \ghv_i(k),
\end{align}
where, $\gtv_i(k)$ is the output of subgradient oracle and,
\begin{equation*}
    \ghv_i(k) = \gtv_i(k) + \sum_{j \in [n] \setminus\{i\}}\hspace{-2mm}\Wv_{ij}(k)\Paren{\Deltav_{ij}(k) + \sum_{l\in [k]}\ntv_{ij}(l)}
\end{equation*}
is the \textit{effective stochastic subgradient} which also includes errors due to quantization and channel noise.

Denote the quantization error at node $i$ by $\qv_i(k) = \sum_{j \in [n]\setminus\{i\}}\Wv_{ij}(k)\Deltav_{ij}(k)$, and noise contribution as $\ntv_i(k) = \sum_{j  \in [n] \setminus \{i\} }\Wv_{ij}(k)\sum_{l\in [k]}\ntv_{ij}(l)$.
From Lemma \ref{lem:unbiased_bounded_variance_quantizer}, 
\begin{align}
    &\Var\Paren{\Br{\qv_i(k)}_s} = \sum_{j  \in [n] \setminus \{i\} }\Var\Paren{\Br{\Wv_{ij}(k)\Deltav_{ij}(k)}_s} \nonumber \\
    &\leq \sum_{j  \in [n] \setminus \{i\} }\beta(k)^2\Pv_{ij}^2\Var\Paren{\Br{\Deltav_{ij}(k)}_s} \leq c_0^2\frac{\Delta^2}{4}.
\end{align}
The last inequality follows from our choice $\beta(k) = c_0k^{-\gamma}$.
For $\alpha(k) = \sqrt{c_1}k^{\tau/2}$, \eqref{eq:accumulated_noise_variance_general} simplifies to $\Var\Paren{\Br{\ntv_i(k)}_s} = $
\begin{equation}
\label{eq:accumulated_noise_variance_geometric_sequence}
    \sigma^2 \frac{c_0^2k^{-2\gamma}}{c_1}\sum_{l\in [k]}l^{-\tau} \leq \sigma^2 \frac{c_0^2k^{-2\gamma}}{c_1}\int_{0}^{k}l^{-\tau}d\tau = \sigma^2\frac{c_0^2}{c_1}\frac{k^{1 - \tau - 2\gamma}}{(1 - \tau)}.
\end{equation}
In particular, \eqref{eq:accumulated_noise_variance_geometric_sequence} implies that the condition of Lemma \ref{lem:convergence_rate_power_requirement_trade-off} is satisfied if the rates for geometric \textit{confidence} and \textit{power control} sequences satisfy $\tau + 2\gamma = 1$.
This gives us $\Var\Paren{\Br{\ntv_i(k)}_s} = \frac{c_0^2\sigma^2}{2\gamma c_1}$.

Since $\ghv_i(k) = \gtv_i(k) + \qv_i(k) + \ntv_i(k)$, taking expectations and recalling Lemma \ref{lem:unbiased_bounded_variance_quantizer}, Asm. \ref{assume:stochastic_subgradient_oracle}, and Asm. \ref{assume:channel_noise_independence}, we have $\Expect\Br{\gvh_i(k)} = \gv_i(k) \in \partial f_i\Paren{\xv_i(k)}$.
We now upper bound the second moments:\footnote{assuming the norm is $\ell_2$, otherwise we need to upper bound the second moment slightly differently.}
\begin{align}
  &\Expect\Br{\norm{\ghv_i(k)}_2^2}
  = \Expect\Br{\norm{\gtv_i(k)}_2^2 + \norm{\qv_i(k)}_2^2 + \norm{\ntv_i(k)}_2^2} \nonumber \\ &+ 2\Expect\Br{\gtv_i(k)^\top\qv_i(k) + \qv_i(k)^\top\ntv_i(k) + \gtv_i(k)^\top\ntv_i(k)}.
\end{align}
From Asm. \ref{assume:channel_noise_independence}, we have
\eas{
\Expect\Br{\qv_i(k)^\top\ntv_i(k)} & = \Br{\Expect\qv_i(k)}^\top\Br{\Expect\ntv_i(k)} = 0 \\
\Expect\Br{\gtv_i(k)^\top\ntv_i(k)} & = \Br{\Expect\gtv_i(k)}^\top\Br{\Expect\ntv_i(k)} = 0.
}
Furthermore, recalling that for a given $i \in [n]$, $\xv_i(k) \in\Fcal_{k-1}$ and $\zv_j(k), \yv_{ij}(k-1) \in \Fcal_{k-1}$ for all $j \in \Ncal(i)$, we have $\Expect\Br{\gtv_i(k)^\top\qv_i(k)} = \Expect\Br{\Expect\Br{\gtv_i(k)^\top\qv_i(k) \vert \Fcal_{k-1}}} = \Expect\Br{\Expect\Br{\gtv_i(k)\vert \Fcal_{k-1}}^\top \Expect\Br{\qv_i(k)\vert \Fcal_{k-1}}} = 0$, where the second equality follows from the fact that given $\Fcal_{k-1}$, $\gtv_i(k)$ and $\qv_i(k)$ are conditionally independent of each other.
Also, 
\ean{
\Expect\Br{\norm{\qv_i(k)}_2^2} & = \Expect\Br{\sum_{s\in [d]}\Br{\qv_i(k)}_s^2} \leq c_0^2 \frac{\Delta^2}{4} d, \\
\Expect\Br{\norm{\ntv(k)}_2^2} & =\Expect\Br{\sum_{s\in [d]}\Br{\ntv_i(k)}_s^2} \leq \frac{c_0^2 \sigma^2 d}{2\gamma c_1}.
}
Accordingly, the effective stochasticity is bounded as:
\begin{equation}
     \Expect\Br{\norm{\gtv_i(k)}_2^2} \leq \Omega^2 + \frac{c_0^2\Delta^2d}{4} + \frac{c_0^2\sigma^2d}{2\gamma c_1} \triangleq \xi^2.
\end{equation}
Now, once the errors due to quantization and channel noise have been clubbed into the stochastic subgradient with bounded effective stochasticity $\xi^2$, analyzing \eqref{eq:dlmd_update_equation_effective_stochasticity} is very similar to distributed dual averaging of \cite{duchi2012TAC}, except for the fact that the consensus matrix $\Wv(k) = \Paren{1 - \beta(k)}\Iv + \beta(k)\Pv$ depends on $k$.
The analysis is \cite{duchi2012TAC} first assumes that there is no stochasticity in the subgradient, i.e. \eqref{eq:dlmd_update_equation_effective_stochasticity} uses $\gv_i(k)$ instead of noisy $\gtv_i(k)$ and derives convergence rate for this case.
It then shows that even if the sub-gradients are noisy, as long as they have bounded variance, much of the analysis remains unchanged. 
We adopt a similar approach.
We first derive an upper bound in Lemma \ref{lem:consensus_error_upper_bound} to the consensus error $\norm{\zvo(k) - \zv_i(k)}_*$ in the presence of  \textit{diminishing confidences} but without any sources of stochasticity, i.e. subgradient oracle is exact, channels are noiseless and there is no data rate constraint.
We then provide an upper bound in Lemma \ref{lem:expected_suboptimality_gap_general} to the expected suboptimality gap $\esf(K)$ in terms of the consensus error $\norm{\zvo(k) - \zv_i(k)}_*$ in the presence of subgradient oracle noise, finite data rate and channel noise. 
Finally to prove Thm. \ref{thm:main_convergence_result}, we combine the results of Lemma \ref{lem:consensus_error_upper_bound} and \ref{lem:expected_suboptimality_gap_general}.

\subsection{Consensus error $\norm{\zvo(k) - \zv_i(k)}_*$ with exact subgradient oracle, infinite data rate and noiseless channel}
\label{subsec:upper_bound_consensus_error}

We first provide an auxillary lemma (Lemma \ref{lem:convergence_doubly_stochastic_matrices}) regarding the rate at which products of doubly stochastic matrices converge to the uniform matrix.
This lemma is required for deriving the upper bound on consensus error in Lemma \ref{lem:consensus_error_upper_bound}.
In other words, if $\Phiv(k,s) = \Wv(k)\Wv(k-1)\ldots\Wv(s)$, then we show that $\Phiv(k,s)$ approaches $\frac{\ones\ones^\top}{n}$ at an exponential rate.

\begin{lemma}
\label{lem:convergence_doubly_stochastic_matrices}
For any fixed doubly stochastic matrix $\Pv \in \Real^{n \times n}$ with eigenvalues $1 = \lambda_1(\Pv) \geq \lambda_2(\Pv) \geq \ldots \geq \lambda_n(\Pv)$, define $\Wv(k)$ as a convex combination of $\Iv$ and $\Pv$: $\Wv(k)  = \Paren{1 - \beta(k)}\Iv + \beta(k) \Pv$, where $\beta(k) = c_0 k^{-\gamma}$ is geometrically decaying confidence sequence, and let $\Phiv(k,s) = \prod_{t \in [s:k]}\Wv(t)$.
Then for any vector $\xv$ in the unit simplex, i.e. $\xv \in \Delta_n = \{\xv \in \Real^n \hspace{1mm} \vert \hspace{1mm} \xv \succcurlyeq 0, \xv^\top\ones = 1\}$, we have:
\begin{align}
     \norm{\Phiv(k,s)\xv - \frac{\ones}{n}}_1 &\leq \sqrt{n}\norm{\Phiv(k,s)\xv - \frac{\ones}{n}}_2 \nonumber\\ 
     &\leq \sqrt{n}e^{-c_0\frac{1-\lambda}{1-\gamma}\Br{\Paren{k+1}^{1-\gamma} - s^{1-\gamma}}},
\end{align}
where, $\lambda = \max\{\lambda_2(\Pv), -\lambda_n(\Pv)\}$ is the eigenvalue of $\Pv$ with second largest magnitude.
\end{lemma}

\begin{IEEEproof}
Note that $\Wv(t), s \leq t \leq k$ are doubly stochastic matrices since convex combination of doubly stochastic matrices remains so.
We have, $\norm{\Phiv(k,s) - \frac{\ones \ones^\top}{n}}_2 = \mu(k,s)$, where $\norm{\cdot}_2$ denotes the operator norm of a matrix, and $\mu(k,s)$ is the second-largest eigenvalue of $\Phi(k,s)$.
From the definition of operator norm, we have for any $\xv \in \Delta_n$, $\norm{\Phiv(k,s)\xv - \frac{\ones}{n}}_2 \leq \mu(k,s)$.
To get an upper bound on $\mu(k,s)$, consider the eigenvalue decomposition of $\Pv = \Vv \Lambdav \Vv^\top$.
Since the eigenspace of $\Iv$ is $\Real^n$, the eigen decomposition of $\Wv(t)$ is $\Wv(t) = \Vv \Paren{\Paren{1 - \beta(t)}\Iv + \beta(t)\Lambdav}\Vv^\top$.
This gives us the eigenvalue decomposition of $\Phiv(k,s)$ to be
\begin{equation*}
    \Phiv(k,s) = \Vv\Br{\prod_{t \in [s:k]}\Paren{\Paren{1 - \beta(t)}\Iv + \beta(t)\Lambdav}}\Vv^\top,
\end{equation*}
The second-largest eigenvalue of $\Phiv(k,s)$ is then $\mu(k,s) = \prod_{t \in [s:k]}\Br{1 + \beta(t)\Paren{\lambda - 1}}$.
Using the fact that $\log\Paren{1 + x} \leq x$,
\begin{equation*}
    \log\mu(k,s) = \hspace{-1.5mm}\sum_{t \in [s:k]}\hspace{-1.5mm}\log\Paren{1 + \beta(t)\Paren{\lambda - 1}} \leq \Paren{\lambda - 1}\sum_{t \in [s:k]}\hspace{-1.5mm}\beta(t),
\end{equation*}
Moreover, $\sum_{t \in [s:k]}\beta(t)$  is evaluated as:
\begin{equation*}
     \sum_{t \in [s:k]}\frac{c_0}{t^{\gamma}} \geq \int_{t=s}^{k+1}\frac{c_0}{t^{\gamma}}dt = \frac{c_0}{1 - \gamma}\Br{\Paren{k+1}^{1- \gamma} - s^{1 - \gamma}}.
\end{equation*}
Substituting this, we get that for $\xv \in \Delta_n$,
\begin{equation*}
    \norm{\Phiv(k,s)\xv - \frac{\ones}{n}}_2 \leq \mu(k,s) \leq e^{-c_0\frac{1 - \lambda}{1 - \gamma}\Br{\Paren{k+1}^{1 - \gamma} - s^{1 - \gamma}}}.
\end{equation*}
This, along with standard norm inequality between $\ell_1$ and $\ell_2$ norms, yields the result.
\end{IEEEproof}

We now present the main lemma of this subsection.
\begin{lemma}
\label{lem:consensus_error_upper_bound}
Consider exact subgradient oracle with noiseless channels having infinite data rate.
Then the state update equations for any node $i \in [n]$, taking into account diminishing confidences is $\zv_i\Paren{k+1} = \sum_{j \in [n] }\Wv_{ij}(k)\zv_j(k) + \gv_i(k)$ and $\xv_i(k+1) = \Pi_{\Xcal}^{\psi}\Paren{\zv_i(k+1), \eta(k)}$.
For a given horizon $K$ and any $k \in [K]$, the network consensus error is upper bounded by
\begin{equation*}
    \norm{\zvo(k) -\zv_i(k)}_* \leq \frac{2L}{c_0\Paren{1 - \lambda}}\Paren{K^{\gamma}\log\Paren{K\sqrt{n}}} + 3L,
\end{equation*}
where $L$ is the Lipschitz constant of each $f_i$.
\end{lemma}
\begin{IEEEproof}
Define $\overline{\Phiv}\Paren{k,s} = \frac{\ones \ones^\top}{n} - \Phiv\Paren{k,s}$. 
Then, unrolling the state update equation from time $k$ back till time $s$, we have
\begin{align}
    \zv_i\Paren{k+1} &=\hspace{-2mm} \sum_{j \in [n] }\Br{\Phiv\Paren{k,s}}_{ji}\zv_j(s) \\ &+ \sum_{r \in [s+1:k]}\Paren{\sum_{j \in [n] }\Br{\Phiv\Paren{k,r}}_{ji}\gv_j\Paren{r-1}} + \gv_i\Paren{k}, \nonumber
\end{align}
where $\Br{\Phiv\Paren{k,s}}_{ji}$ denotes the $j^{th}$ entry of the $i^{th}$ column of $\Phiv\Paren{k,s}$. The average state evolves as:
\begin{align}
    &\zvo\Paren{k+1} = \frac{1}{n}\sum_{i \in [n] }\zv_i\Paren{k+1} \\ &\hspace{-2mm}= \frac{1}{n}\sum_{i \in [n] }\hspace{-1mm}\Paren{\sum_{j \in [n] }\hspace{-1mm}\Wv_{ij}(k)\zv_j(k)+\gv_i\Paren{k}} = \zvo\Paren{k} +\frac{1}{n}\hspace{-0.5mm}\sum_{i \in [n] }\gv_i\Paren{k}. \nonumber
\end{align}
Since we initialize our dual iterates to $\zv_i\Paren{0} = 0$ for all $i \in [n]$, unrolling the above yields $\zvo\Paren{k} - \zv_i\Paren{k} = $
\begin{align}
\label{eq:consensus_error_expression_1}
     &\sum_{s \in [k-1]}\sum_{j \in [n] }\Paren{\frac{1}{n} - \Br{\Phiv\Paren{k-1,s}}_{ji}}\gv_j\Paren{s-1} \nonumber\\ &+ \Paren{\frac{1}{n}\sum_{j \in [n] }\Paren{\gv_j\Paren{k-1} - \gv_i\Paren{k-1}}}
\end{align}

Furthermore, taking $\norm{\cdot}_*$ on both sides of \eqref{eq:consensus_error_expression_1}, and using the fact that each component function $f_i$ of our objective is $L$-Lipschitz continuous, i.e. $\norm{\gv_i\Paren{k}}_* \leq L$, simple algebraic manipulation gives us:
\begin{equation}
\label{eq:consensus_error_sum}
    \norm{\zvo(k) - \zv_i(k)}_* \leq \sum_{s \in [k-1]}L\norm{\Phiv(k-1,s)\ev_i - \frac{\ones}{n}}_1 + 2L.
\end{equation}
Here, $\ev_i$ denotes the $i^{th}$ standard basis vector.
We now split the sum in \eqref{eq:consensus_error_sum} into two parts separated by a point $\kh$.
The first part consists of iteration steps $k$ for which $\norm{\Phiv(k-1,s) - \frac{\ones}{n}}_1$ is small, and the second part has limited number of terms in the summation.
Since $\ev_i \in \Delta_n$, from Lemma \ref{lem:convergence_doubly_stochastic_matrices} we have, $\norm{\Phiv\Paren{k-1,s}\ev_i - \frac{\ones}{n}}_1 \leq \sqrt{n}e^{-c_0\frac{1-\lambda}{1 - \gamma}\Br{k^{1-\gamma} - s^{1-\gamma}}}$.
If we split the sum \eqref{eq:consensus_error_sum} at $\kh = \Br{k^{1-\gamma} - \frac{1 - \gamma}{c_0\Paren{1 -\ \lambda}}\loge{K\sqrt{n}}}^{\frac{1}{1-\gamma}}$, we have for $s \leq \kh$, $\norm{\Phiv(k-1,s)\ev_i - \frac{\ones}{n}}_1 \leq \frac{1}{K}$, and for larger $s$, we use the more trivial bound $\norm{\Phiv(k-1,s)\ev_i - \frac{\ones}{n}}_1 \leq 2$. 
Doing this yields,
\begin{equation}
    \label{eq:consensus_error_sum_2}
    \norm{\zvo(k) - \zv_i(k)}_* \leq L\sum_{s = \kh + 1}^{k-1}\norm{\Phiv\Paren{k-1,s}\ev_i - \frac{\ones}{n}}_1 + 3L.
\end{equation}
Using the inequality $\Paren{1-x}^n \geq 1 - nx$ when $x \leq 1$ for any even positive integer $n$, we get,
\begin{equation*}
    \Paren{k-1} - \Paren{\kh - 1} + 1 \leq \frac{k^{\gamma}\loge{K\sqrt{n}}}{c_0\Paren{1 - \lambda}}
\end{equation*}
Substituting this in \eqref{eq:consensus_error_sum_2} and recalling that $k \leq K$, we obtain the desired result.
\end{IEEEproof}

\subsection{Expected suboptimality gap in the presence of stochastic subgradient oracle and noisy channel with finite data rate}
\label{subsec:expected suboptimality_gap_general}

\begin{lemma}
\label{lem:expected_suboptimality_gap_general}
Let the sequences $\{\xv_i(k)\}_{k\in \Nbb}$ and $\{\zv_i(k)\}_{k  \in \Nbb}$ be generated by \gls{dlmd}.
Then, for any $\xv^* \in \Xcal$ and for each node $i \in [n]$, we have $\Expect f(\xvh_i(K)) - f(x^*) \leq$
\begin{align*}
    &\frac{1}{K\eta(K)}\psi(\xv^*) + \frac{\xi^2}{2K}\sum_{k \in [K]}\eta(k-1) \\ &+ \frac{L + \xi}{nK}\sum_{k \in [K]}
    \sum_{i \in [n] }\eta(k)\Expect\norm{\zvo(k) - \zv_i(k)} \\ &+ \frac{L}{K}\sum_{k \in [K]}\eta(k)\Expect\norm{\zvo(k) - \zv_i(k)},
\end{align*}
where $\norm{\cdot}$ denotes the $\ell_2$-norm which is its own dual.
\end{lemma}
\begin{IEEEproof}
The proof of this follows the proof of \cite[Thm. 1]{duchi2012TAC} very closely.
Hence, we highlight the major differences from that proof and omit the tedious algebraic manipulations while referring the reader to the relevant equations in \cite{duchi2012TAC} for reference. 
Consider the running sum $S(K) = \sum_{k \in [K]}\Br{f(\yv(k)) - f(\xv^*)}$.
Following steps similar to \cite[Eqns. 19, 20]{duchi2012TAC},
\begin{align}
    \label{eq:lemma8_proof_step1}
 S(K) &\leq 
 \frac{1}{n}\sum_{k \in [K]}\sum_{i \in [n] }\inprod{\ghv_i(k), \xv_i(k) - \xv^*} \nonumber \\
 & + \frac{L}{n}\sum_{k \in [K]}\sum_{i \in [n] }\norm{\yv(k) - \xv_i(k)} \nonumber\\ 
 & + \frac{1}{n}\sum_{k \in [K]}\sum_{i \in [n] }\inprod{\gv_i(k) - \ghv_i(k), \xv_i(k) - \xv^*}.
\end{align}
The expression inside the summation over $k$ in the first term of \eqref{eq:lemma8_proof_step1} can be written as $\inprod{\frac{1}{n}\sum_{i \in [n] }\ghv_i(k), \xv_i(k) - \xv^*} =$
\begin{equation}
    \label{eq:lemma8_proof_step2} 
\inprod{\frac{1}{n}\sum_{i \in [n] }\ghv_i(k), \yv(k) - \xv^*} 
+\inprod{\frac{1}{n}\sum_{i \in [n] }\ghv_i(k), \xv_i(k) - \yv(k)}.
\end{equation}

Note that even in the presence of stochasticity, the centralized state sequence $\zvo(k)$ evolves in a very simple way:
\begin{align*}
    \zvo(k+1) &= \frac{1}{n}\sum_{i \in [n] }\zv_i(k+1) \\
    &= \frac{1}{n}\sum_{i \in [n] }\Paren{\sum_{j \in [n] }\Wv_{ij}(k)\zv_j(k) + \ghv_i(k)} \\
    &= \zvo(k) + \frac{1}{n}\sum_{i \in [n] }\ghv_i(k).
\end{align*}

From Lemma \ref{lem:projection_sequence_lemma}, with the sequence $\left\{\frac{1}{n}\sum_{i \in [n]}\ghv_i(k)\right\}$ in place of $\{\gv_i(k)\}$ we have $\sum_{k \in [K]}\inprod{\frac{1}{n}\sum_{i \in [n] }\ghv_i(k), \yv(k) - \xv^*}$
\begin{equation}
    \leq \frac{1}{\eta(K)}\psi(\xv^*) + \f 12\sum_{k \in [K]}\eta(k-1)\norm{\frac{1}{n}\sum_{i \in [n] }\ghv_i(k)}^2 ,
\end{equation}
Taking expectations and using Holder's inequality,
\begin{align}
    &\Expect\Br{\norm{\frac{1}{n}\sum_{i \in [n] }\ghv_i(k)}^2} \leq \frac{1}{n^2}\sum_{i,j \in [n]}\Expect\Br{\norm{\ghv_i(k)}\norm{\ghv_j(k)}} \nonumber \\
    &\leq \frac{1}{n^2}\sum_{i,j \in [n]}\Paren{\Expect\Br{\norm{\ghv_i(k)}^2}}^{\f 12}\Paren{\Expect\Br{\norm{\ghv_j(k)}^2}}^{\f 12} \leq \xi^2.
\end{align}
The second term in \eqref{eq:lemma8_proof_step2} is, $\Expect\Br{\inprod{\ghv_i(k),\xv_i(k) - \yv(k)}}$
\begin{align}
      &\leq \Expect\Br{\norm{\ghv_i(k)}\norm{\xv_i(k) - \yv(k)}} \nonumber \\
      &= \Expect\Br{\Expect\Br{\norm{\ghv_i(k)}\vert \Fcal_{k-1}}\cdot \norm{\xv_i(k) - \yv(k)}} \nonumber \\
      &\leq \xi \Expect\Br{\norm{\xv_i(k) - \yv(k)}}.
\end{align}
We then get, $\sum_{k \in [K]}\frac{1}{n}\sum_{i \in [n] }\Expect\Br{\inprod{\ghv_i(k), \xv_i(k) - \xv^*}}$
\begin{align}
     &\leq \frac{1}{\eta(K)}\psi(\xv^*) + \frac{\xi^2}{2}\sum_{k \in [K]}\eta(k-1) \nonumber \\
     &+ \frac{\xi}{n}\sum_{k \in [K]}\sum_{i \in [n] }\eta(k)\Expect\Br{\norm{\zvo(k) - \zv_i(k)}},
\end{align}
where the final term is obtained by taking expectations on both sides of Lemma \ref{lem:projection_lipschitz_continuity}.
The second term of \eqref{eq:lemma8_proof_step1} is upper bounded similarly using Lemma \ref{lem:projection_lipschitz_continuity}, that is, $\frac{L}{n}\sum_{k \in [K]}\sum_{i \in [n] }\Expect\Br{\norm{\yv(k) - \xv_i(k)}} \leq \frac{L}{n}\sum_{k \in [K]}\sum_{i \in [n] }\eta(k)\Expect\Br{\norm{\zvo(k) - \zv_i(k)}}$.

Finally, taking expectations of the third term in \eqref{eq:lemma8_proof_step1}, $\Expect\Br{\inprod{\gv_i(k) - \gvh_i(k), \xv_i(k) - \xv^*}}$
\begin{align*}
    &= \Expect\Br{\Expect\Br{\inprod{\gv_i(k) - \gvh_i(k), \xv_i(k) - \xv^*}\vert \Fcal_{k-1}}}\\
    &= \Expect\Br{\inprod{\Expect\Br{\gv_i(k) - \gvh_i(k) \vert \Fcal_{k-1}}, \xv_i(k) - \xv^*}} = 0
\end{align*}
With all this, $\Expect \Br{S(K)}$ eventually boils down to
\begin{align}
    \label{eq:lemma8_proof_step3}
    \Expect \Br{S(K)} &\leq \frac{1}{\eta(K)}\psi(\xv^*) + \frac{\xi^2}{2}\sum_{k \in [K]}\eta(k-1) \nonumber \\ &+ \frac{L + \xi}{n}\sum_{k \in [K]}\sum_{i \in [n] }\eta(k)\Expect\norm{\zvo(k) - \zv_i(k)}.
\end{align}
Since $\xvh = \frac{1}{K}\sum_{k \in [K]}\xv_i(k)$, using Jensen's inequality,
\begin{equation*}
    \Expect\Br{f(\xh_i(K))} - f(\xv^*) \leq \frac{1}{K}\Expect\Br{\sum_{k \in [K]}f(\xv_i(k)) - f(\xv^*)}
\end{equation*}
Finally, recalling Lemma \ref{lem:centralized_individual_sequence_relation}, 
$$\Expect f_i(\xvh_i(k)) - f(\xv^*) \leq \Expect\Br{S(K)} + L\hspace{-2mm}\sum_{k \in [K]}\eta(k)\norm{\zvo(k) - \zv_i(k)}.$$
Substituting \eqref{eq:lemma8_proof_step3} here, the proof is complete.
\end{IEEEproof}

\subsection{Proof of Theorem \ref{thm:main_convergence_result}}

The proof of Thm. \ref{thm:main_convergence_result} is essentially substituting Lemma \ref{lem:consensus_error_upper_bound} in Lemma \ref{lem:expected_suboptimality_gap_general}.
First, note that $\Expect\Br{\norm{\ghv_i(k)}} \leq \xi$.
This follows from the fact that $\norm{\cdot}$ is a convex function and Jensen's inequality yields $\Paren{\Expect\Br{\norm{\ghv_i(k)}}}^2 \leq \Expect\Br{\norm{\ghv_i(k)}^2} \leq \xi^2$.
This shows that an analogous version of Lemma \ref{lem:consensus_error_upper_bound} holds even in the presence of stochasticity.
Simply replacing $\norm{\zvo(k) - \zv_i(k)}$ by $\Expect\norm{\zvo(k) - \zv_i(k)}$, and $\norm{\gv_i(k)} \leq L$ everywhere in the derivation with $\Expect\norm{\ghv_i(k)} \leq \xi$, we get
\begin{equation}
    \label{eq:expected_consensus_error_upper_bound}
    \Expect\norm{\zvo(k) - \zv_i(k)} \leq \frac{2 \xi}{c_0(1 - \lambda)}K^{\gamma}\loge{K\sqrt{n}} + 3\xi.
\end{equation}
Replacing the expected network consensus error $\Expect\norm{\zvo(k) - \zv_i(k)}$ in the expression of Lemma \ref{lem:expected_suboptimality_gap_general}  by this upper bound in \eqref{eq:expected_consensus_error_upper_bound}, the expression for stepsize $\eta(k)$, the inequality $\sum_{k \in [K]}k^{-p} \leq \int_{0}^{K}k^{-p}dk = K^{1-p}/(1-p)$, and the facts that $1 - \lambda < 1$ and $L \leq \xi$, we get our result.

\section{Proof of Theorem \ref{thm:quantizer_unsaturation}}
\label{appendix:proof_of_theorem_2}

Let us denote the input to the quantizer along the link between nodes $i$ and $j$ by
\begin{equation}
    \omegav_{ij}(k) = \zv_j(k) - \yv_{ij}(k-1).
\end{equation}
We need to ensure that $\omegav_{ij}(k)$ remains within the dynamic range of the quantizer with high-probability.
Note that:
\begin{align*}
    &\omegav_{ij}(k) = \Br{\zv_j(k-1) - \yv_{ij}(k-1)} \\
    &- \beta(k-1)\hspace{-2mm}\sum_{w \in \Ncal(j)}\hspace{-2mm}\Pv_{jw}\Paren{\zv_j(k-1)- \ytv_{jw}(k-1)} + \gtv_j(k-1).
\end{align*}
Some simple algebraic manipulations yield:
\begin{align}
    \omegav_{ij}(k) &= -\Deltav_{ij}(k-1) \nonumber \\
    &+ \beta(k-1)\sum_{w \in \Ncal(j)}\hspace{-2mm}\Pv_{jw}\hspace{-1mm}\Br{\sum_{s \in [k-1]}\hspace{-2mm}\ntv_{jw}(s) + \Deltav_{jw}(k-1)} \nonumber \\ &+ \gtv_j(k-1).
\end{align}
Recall that conditioned on the event that all quantizers were unsaturated at time $k-1$, quantization, channel noise, and stochastic subgradient oracle noise are unbiased.
So, we have:
\begin{align*}
    &\Expect\Br{\omegav_{ij}(k)} = \gv_j(k-1) \nonumber\\ \implies &\Expect\Br{\omegav_{ij}}_p = \Br{\gv_j(k-1)}_p \leq L, \hspace{2mm} \text{for} \hspace{2mm} p \in [d].
\end{align*}
From the above expression, making use of double stochasticity of $\Pv$ and plugging in the geometric sequences for $\beta(k)$ and $\alpha(k)$, we also get $\Var\Br{\omegav_{ij}(k)}_p$
\begin{align}
     &\leq  \Paren{1+\beta^2(k-1)}\frac{\Delta^2}{4} + \frac{\beta^2(k-1)\sigma^2}{c_1}\sum_{s \in [k-1]}\hspace{-2mm} s^{-\tau} + \Omega^2 \nonumber \\
     &\leq \Paren{\frac{c_0^2\Delta^2}{2} + \frac{c_0^2\sigma^2}{2\gamma c_1} + \Omega^2} \triangleq \Psi^2.
\end{align}
where the last inequality makes a simplifying assumption $c_0 \geq 1$.
Note that in the absence of any stochasticity, a quantizer range of $L$ would have ensured that it never gets saturated.
But in the presence of quantization, channel noise and stochasticity from subgradient oracle, no matter what the chosen range of the quantizer is, there will always be a non-zero probability that it gets saturated.
However, if we choose the dynamic range of the quantizers to be large enough, as long as it is at least $U_{\min} = L$, we can upper bound the probability of quantizer saturation by an arbitrarily small value.
From Chebyshev's inequality, for any $p \in [d]$, we have:
\begin{align}
    &\Prob\Paren{|\Br{\omegav_{ij}(k)}_p - \Expect \Br{\omegav_{ij}(k)}_p| \geq u} \leq \frac{\Psi^2}{u^2} \nonumber\\ \implies &\Prob\Paren{-U \leq \Br{\omegav_{ij}(k)}_p \leq +U} \geq 1 - \frac{\Psi^2}{\Paren{U - U_{\min}}^2},
\end{align}
where, $U$ is the chosen dynamic range of the quantizer, $U_{\min} = L$ is the Lipschitz constant of the objective function, and $\Psi^2$ is the variance of each coordinate of the vector input to the quantizer.
Using the chain rule of probability, $\Prob\Paren{\norm{\omegav_{ij}(k)}_{\infty} \leq U \hspace{1mm} \forall \hspace{1mm} i,j} = $
\begin{small}
\begin{align}
    &\Prob\Paren{\norm{\omegav_{ij}(k)}_{\infty} \leq U \hspace{1mm} \forall \hspace{1mm} i,j \hspace{1mm} \vert \hspace{1mm} \norm{\omegav_{ij}(k-1)}_{\infty} \leq U \hspace{1mm} \forall \hspace{1mm} i,j} \nonumber\\ \cdot &\Prob\Paren{\norm{\omegav_{ij}(k-1)}_{\infty} \leq U \hspace{1mm} \forall \hspace{1mm} i,j} \nonumber \\
    = &\prod_{s\in [k]}\Prob\Paren{\norm{\omegav_{ij}(s)}_{\infty} \leq U \hspace{1mm} \forall \hspace{1mm} i,j \hspace{1mm} \vert \hspace{1mm}\norm{\omegav_{ij}(s-1)}_{\infty} \leq U \hspace{1mm} \forall \hspace{1mm} i,j} \nonumber \\
    \geq &\Paren{1 - \frac{\Psi^2}{\Paren{U - U_{\min}}^2}}^{2kd \vert E \vert}.
\end{align}
\end{small}
where, the $2$ appears in the exponent because the network links are taken to be bidirectional.
This proves Thm. \ref{thm:quantizer_unsaturation}.

\end{document}